\documentclass[11pt]{amsart}
\usepackage{amssymb}
\usepackage{amsmath}
\usepackage{mathabx}
\usepackage[latin1]{inputenc}
\usepackage{enumerate}
\usepackage{mathtools}
\usepackage{color}
\usepackage{hyperref}

\def\a{\alpha}

\def\e{\varepsilon}

\let\newpf\proof \let\proof\relax 
\newenvironment{pf}{\newpf[\proofname]}{\qed\endtrivlist}

\newcommand{\ba}{\overline{A}}

\def\be{\begin{equation}}
\def\ee{\end{equation}}

\def\ba{{\begin{align}}}
\def\ea{{\end{align}}}

\def\bm{\begin{matrix}}
\def\em{\end{matrix}}

\renewcommand{\sl}{{\mathrm{sl}}}

\def\u{{\mathbb U}}

\def\a{{\alpha}}

\def\0{{\mathbf 0}}

\newtheorem{Theorem}{Theorem}[section]
\newtheorem{Lemma}{Lemma}[section]
\newtheorem{Proposition}{Proposition}[section]
\newtheorem{Corollary}{Corollary}[section]
\newtheorem{Remark}{Remark}[section]

\newtheorem{Definition}{Definition}[section]

\numberwithin{equation}{section}

\theoremstyle{definition}

\newcommand{\id}{\operatorname{id}}

\def\tF{\widetilde{F}}

\newcommand{\C}{{\mathbb C}}

\newcommand{\N}{{\mathbb N}}
\newcommand{\Q}{{\mathbb Q}}
\newcommand{\R}{{\mathbb R}}
\newcommand{\T}{{\mathbb T}}

\newcommand{\Z}{{\mathbb Z}}

\newcommand{\la}{\langle}
\newcommand{\ra}{\rangle}

\def\B0{{\bold{0}}}

\catcode`\@=12

\def\Empty{}
\newcommand\oplabel[1]{
  \def\OpArg{#1} \ifx \OpArg\Empty {} \else
    \label{#1}
  \fi}

\newcommand{\comm}[1]{}
\newcommand{\comment}[1]{}

\begin{document}

\title[]{Exponential Dynamical Localization: Criterion and  Applications}

\author{Linrui Ge}
\address{
Department of Mathematics, Nanjing University, Nanjing 210093, China
}

 \email{lingruige10@163.com}

\author {Jiangong You}
\address{
Chern Institute of Mathematics and LPMC, Nankai University, Tianjin 300071, China} \email{jyou@nju.edu.cn}

\author{Qi Zhou}
\address{
Department of Mathematics, Nanjing University, Nanjing 210093, China
}

 \email{qizhou628@gmail.com, qizhou@nju.edu.cn}

 \begin{abstract}
 We give a criterion for exponential dynamical localization in expectation (EDL) for ergodic families of operators acting on $\ell^2(\Z^d)$. As applications, we prove EDL for a class of quasi-periodic long-range operators on $\ell^2(\Z^d)$.
\end{abstract}

\maketitle

\section{Introduction}

%
%
%
%
%

Localization of particles and waves in disordered media is one of the most intriguing phenomena in solid-state physics. Date back to 1958, Anderson \cite{anderson} firstly used a tight-binding model of an electron in a disordered lattice to explain the physical phenomenon that anomalously long relaxation times of electron spins in doped semiconductors. Since that ground breaking work,  physicists found Anderson's discovery  plays an important  role in the quantization of Hall conductance \cite{ag,h,ntw},  in the emerging subject of optical crystals \cite{fk}, e.t.c. One may  consult \cite{history} for more about the  history of the localization theory.

The mathematical models of the above problems appear often  as ergodic families of Schr\"odinger operators. Let $(\Omega,d\mu, S)$ be a measurable space with a Borel probability measure, equipped with an ergodic family of maps $S=\{S_n,n\in\Z^d\}$ from  $\Omega$ to $\Omega$ such that $S_{m+n}=S_mS_n$.  A measurable family $(H_\omega)_{\omega\in\Omega}$ of bounded linear self-adjoint operators on $\ell^2(\Z^d)$ is called ergodic if $H_{S_n\omega}=T_nH_\omega T_{-n}$ where $T_n$ is the translation in $\Z^d$ by the vector $n$. Motivated by physical backgrounds, in the past sixty years, localization property of ergodic families of operators $(H_{\omega})_{\omega\in \Omega}$   has been extensively studied.
There are several different  mathematical definitions for localization from weak sense to strong sense. The weakest one is Anderson localization (AL):  $H_{\omega}$ is said to display Anderson localization for $a.e.$ $\omega$, if $H_{\omega}$ has pure point spectrum with exponentially decaying eigenfunctions for $a.e.$ $\omega$. However, in physics, localization often means dynamical localization (DL), i.e.,  the wave-packets under the Schr\"odinger time evolution $e^{itH_\omega}$ keep localized if the initial wave-packet is localized.   One generally accepted definition  of DL is the following:
  for any $q>0$,
\begin{equation}\label{sdl single}
\sup\limits_{t}\sum\limits_{n\in\Z}(1+|n|)^q|\langle\delta_n,e^{-itH_{\omega}}\delta_0\rangle|<C(\omega)<\infty.
\end{equation}

When considering an ergodic family operators $(H_{\omega})_{\omega\in\Omega}$, the above defined DL often holds only for a full measure of $\omega$, moreover $C(\omega)$ is not uniform in $\omega$.  So it is natural to consider the dynamical localization  in expectation. The
 strongest dynamical localization is  the \textit{exponential localization in expectation} (EDL)  \cite{ag,jk,jkl}:
\begin{equation}\label{edl}
\int_{\Omega}\sup\limits_{t\in\R}|\langle\delta_k, e^{-itH_\omega}\delta_{\ell}\rangle|d\mu \leq Ce^{-\gamma|k-\ell|},
\end{equation}
where
$$\gamma(H):= \liminf_{n\rightarrow \infty} (-\frac{1}{|n|} \ln \int_{\Omega}\sup\limits_{t\in\R}|\langle\delta_n, e^{-itH_\omega}\delta_{0}\rangle|d\mu)  $$  is defined  as the \textit{exponential decay rate in expectation} \cite{jk,jkl}. This quantity is useful since it is  obviously
connected to the minimal inverse correlation length.
  As pointed out in \cite{ag,jk}, EDL leads to various interesting physical conclusions, for example, the exponential decay of the two-point function at the ground state and positive temperatures with correlation length staying uniformly bounded as temperature goes to zero. Moreover, it is indeed EDL that is often implicitly assumed as manifesting localization in physics literature. That makes it particularly interesting to establish EDL for physically relevant models.

Dynamical localization  has already been well studied for the random  cases \cite{a,am,asfh,dks,ks}. So we will focus our attention on  the quasi-periodic models, which even have stronger  backgrounds in physics \cite{aos,oa}.  As pointed out in \cite{gj}, establishing dynamical localization in the quasi-periodic case requires approaches that are quite different from that of the random case (which usually use Aizenman-Molchanov's fractional moments method \cite{am}).
So far, there are few EDL  results for one dimensional quasi-periodic Schr\"odinger operators and there is no EDL result for higher dimensional quasi-periodic Schr\"odinger operators.  In this paper, we will provide a general criterion for EDL by information of eigenfunctions and apply it to some popular quasi-periodic models.  More applications to other quasiperiodic models will be given in our forthcoming paper. Hopefully, the criterion can also be used to deal with other ergodic Schr\"odinger operators.

Date back to 1980s, the most studied quasiperiodic models have the following form:
\begin{equation}\label{longmoregeneral}
(L_{V, \lambda W, \alpha,\theta}u)_n=\sum\limits_{k\in\Z^d}V_ku_{n-k}+ \lambda W (\theta+ \langle n,\alpha\rangle)u_n,
\end{equation}
where  $\alpha\in\T^d$ is the frequency, $V_k\in\C$ is the Fourier coefficient of  a real analytic function $V:\T^d \rightarrow  \R$, $W$ is a continuous function defined on $\T$. It is an ergodic  self-adjoint operator which is defined on $\ell^2(\Z^d)$. There are two fundamental results when considering the localization property of $L_{V,  \lambda W, \alpha,\theta}$.  Inspired by  pioneer work  of  Fr\"ohlich-Spencer-Wittwer \cite{fsw} and Sinai \cite{Sinai},  Chulaevsky-Dinaburg \cite{cd} proved that if $W$ is a cosine-like function, then for any fixed phase $\theta$, and for positive measure $\alpha$, $L_{V,  \lambda W, \alpha,\theta}$ has AL for sufficiently large coupling constant $\lambda$. Bourgain \cite{b1} generalized this result to arbitrary real analytic function $W$.

 If $W(\theta)=2\cos(2\pi \theta)$, then \eqref{longmoregeneral} reduces to the famous quasiperiodic long-range operators:
\begin{equation}\label{longgeneral}
(L_{V, \lambda, \alpha,\theta}u)_n=\sum\limits_{k\in\Z^d}V_ku_{n-k}+ 2\lambda \cos (\theta+ \langle n,\alpha\rangle)u_n.
\end{equation}
The operator (\ref{longgeneral}) has received a lot of attention \cite{aj1,ayz,ayz1,bj02} since it   is  the Aubry duality of the general one dimensional quasiperiodic Schr\"odinger operator
\begin{equation}\label{1.5}(H_{\lambda^{-1} V,\alpha, \theta}  u)_n= u_{n+1}+u_{n-1} + \lambda^{-1} V(\theta_1+n\alpha_1,\cdots,\theta_d+n\alpha_d)u_n,\end{equation} and thus  (\ref{longgeneral}) carries enormous information of (\ref{1.5})
\cite{cd1,gjls}.

 If furthermore $V(\theta)=2\cos(2\pi \theta)$, then \eqref{longgeneral} reduces to the most famous almost Mathieu operator (AMO):
\begin{equation}\label{amo}
(H_{\lambda, \alpha,\theta}u)_n= u_{n+1}+ u_{n-1} +  2\lambda\cos(\theta+ n\alpha) u_n,
\end{equation}
 The AMO was first introduced by Peierls \cite{p},
as a model for an electron on a 2D lattice, acted on by a homogeneous
magnetic field \cite{harper,R}. This model
has been extensively studied not only because of  its importance
in  physics \cite{aos,oa},
but also as a fascinating mathematical object. More detailed property of AMO will be discussed separately in Subsection \ref{sec-amo}.

\subsection{EDL for quasi-periodic long-range operators }
Considering \eqref{longgeneral}, if $d=1$,  Bourgain-Jitomirskaya \cite{bj02} proved that for  any fixed Diophantine number $\alpha$\footnote{ $\alpha \in\T^d$ is  called {\it Diophantine}, denoted by $\alpha \in {\rm DC}_d(\kappa',\tau)$, if there exist $\kappa'>0$ and $\tau>d-1$ such that
\begin{equation}\label{dio}
{\rm DC}_d(\kappa',\tau):=\left\{\alpha \in\T^d:  \inf_{j \in \Z}\left| \langle n,\alpha  \rangle - j \right|
> \frac{\kappa'}{|n|^{\tau}},\quad \forall \  n\in\Z^d\backslash\{0\} \right\}.
\end{equation}
Let ${\rm DC}_d:=\bigcup_{\kappa'>0,
\, \tau>d-1} {\rm DC}_d(\kappa',\tau)$.}, $L_{V,\lambda,\alpha,\theta}$ has DL for sufficiently large $\lambda$ and a.e. $\theta$. This result is non-perturbative in the sense that the largeness of $\lambda$ doesn't  depend on the Diophantine constant of $\alpha$.
%
%
%

If $d\geq 2$,  Jitomirskaya and Kachkovskiy   \cite{jk} proved that  for fixed $\alpha\in DC_d$,   $L_{ V,\lambda,\alpha,\theta}$ has pure point spectrum for large enough $\lambda$ and $a.e.$ $\theta$.   Compared to the one dimensional result, here one can only obtain the perturbative  result, i.e., the largeness of $\lambda$ depends on $\alpha$.

We remark that Jitomirskaya-Kachkovskiy's result \cite{jk} are based on Eliasson's full measure reducibility results \cite{E92} where the cocycles are full measure reducible only in $C^\omega$ topology.  By Aubry duality,  although all the eigenfunctions decay exponentially,  they don't have uniform decay. That's the reason why  
they could only obtain pure point spectrum instead of AL.  What we will prove is,  under the exact same setting as in \cite{jk},  the family of operator $(L_{V,\lambda,\alpha,\theta})_{\theta\in\T}$ not only displays AL but also EDL. Now we precisely formulate our result.

Recall that for a bounded analytic (possibly matrix valued) function $F$ defined on $ \{ \theta |  | \Im \theta |< h \}$, let
$
\lvert F\rvert _h=  \sup_{ | \Im \theta |< h } \| F(\theta)\| $, and denote by $C^\omega_{h}(\T^d,*)$ the
set of all these $*$-valued functions ($*$ will usually denote $\R$, $sl(2,\R)$
$SL(2,\R)$).  Denote $C^{\omega}(\T^d,\R)$ by the union $\cup_{h>0}C_h^{\omega}(\T^d,\R)$.

\begin{Theorem}\label{long}
If $\alpha\in DC_d$,  $h>0$, $V\in C_h^{\omega}(\T^d,\R)$. Then for any $\e>0$, there exists $\lambda_0(\alpha,V,d,\e)$, such that
$L_{V,\lambda,\alpha,\theta}$ has EDL with $\gamma(L) \geq 2\pi(h-\e)$ if $\lambda>\lambda_0$.
\end{Theorem}

\begin{Remark}
To the best knowledge of the authors, Theorem \ref{long} gives the first dynamical localization type result for multidimensional quasi-periodic  operators.
\end{Remark}

\begin{Remark}
If $d=1$, the result we get is non-perturbative, see Corollary \ref{nonperturbativelong}.
\end{Remark}

\subsection{EDL for quasi-periodic Schr\"odinger  operators on $\ell^2(\Z^d)$}
If we take $V=\sum_{i=1}^d2\cos2\pi{x_i}$,  \eqref{longmoregeneral} is reduced to  quasi-periodic Schr\"odinger operators on $\ell^2(\Z^d)$:
\begin{equation}\label{highschrodinger-1}
H_{\lambda W,\alpha,\theta}= \Delta +  \lambda W(\theta+\langle n,\alpha\rangle)\delta_{nn'},
\end{equation}
where $\Delta$ is the usual Laplacian on $\Z^d$ lattice.

If $d=1$,  Eliasson \cite{Eli97} proved that for any {\it fixed} Diophantine frequency, $H_{\lambda W,\alpha,\theta}$ has pure point for $a.e.$ $\theta$ and large enough $\lambda.$   Bourgain and Goldstein \cite{bg} proved that in the positive Lyapunov exponent regime, for any {\it fixed} phase, $H_{\lambda W,\alpha,\theta}$ has AL for $a.e.$ Diophantine frequency.   Klein \cite{K} generalized the results in \cite{bg} to more general Gevrey potentials.

If $d=2$, Bourgain, Goldstein and Schlag \cite{bgs2} proved that  for any {\it fixed} $\theta$, $H_{\lambda W,\alpha,\theta}$ has AL for sufficiently large $\lambda$  and a positive Lebesgue measure set of $\alpha$.   As a matter of  fact, Bougain-Godstein-Schlag's method \cite{bgs2}  works for more complicated potential $W$ of the form $$W(n_1,n_2)=w(\theta_1+n_1\alpha_1, \theta_2+n_2\alpha_2)$$ where $w$ is a real analytic function on $\T^2$  which is non-constant on any horizontal or vertical line. Later, Bourgain \cite{b} generalize the result \cite{bgs2} to the case  $d\geq 3$.

However, even the  \textit{weakest} dynamical localization result of (\ref{highschrodinger-1}) is not known. In this paper, we will establish the first  EDL result for a family of multidimentional quasi-periodic  Schr\"odinger  operator: 

\begin{Theorem}\label{highamo}
	For $\alpha\in {\rm DC}_d$ and any $\e>0$, there exists $\lambda_0(\alpha,d,\e)$, such that if $\lambda>\lambda_0$, then \begin{equation}\label{highschrodinger}
	L_{\lambda, \alpha,\theta}= \Delta + 2\lambda \cos(\theta+\langle n,\alpha\rangle)\delta_{nn'},
	\end{equation} has EDL with exponential decay rate in expectation $\gamma(L)\geq (1-\e)\ln\lambda$.
\end{Theorem}

\subsection{EDL for Almost Mathieu Operators} \label{sec-amo}

Theorem \ref{highamo} established EDL for $L_{\lambda,\alpha,\theta}$ with large coupling constant $\lambda$. When $d=1$, (\ref{highschrodinger}) reduces the famous almost Mathieu operators $H_{\lambda,\alpha,\theta}$. In this case EDL can be proved for  $|\lambda|>1$:

\begin{Theorem}\label{amo}
Let $\alpha\in DC_1$,  $|\lambda|>1$, then almost Mathieu operators $H_{\lambda,\alpha,\theta}$ has EDL with  exponential decay rate in expectation  $\gamma(H)=\ln|\lambda|$.
\end{Theorem}

We briefly review the localization type results for AMO, if $\alpha\in DC_1$,
Jitomirskaya \cite{J} proved that if $\lambda>1$  then $H_{\lambda,\alpha,\theta}$  has AL for $\theta\in \Theta$. Jitomirskaya and Liu \cite{JLiu1} proved sharp transition in the arithmetics of phase between localization and singular continuous spectrum. What's remarkable is that  they uncover a reflective-hierarchy structure for quasi-periodic eigenfunctions.
If $\beta(\alpha)>0$, Avila and Jitomirskaya \cite{aj} proved that for $\lambda>e^{\frac{16\beta}{9}}$, $H_{\lambda,\alpha,\theta}$ has AL for $\theta\in\Theta$. Here,  $\beta(\alpha)$ is defined as
\begin{equation*}\beta(\alpha):=\limsup_{n\rightarrow \infty}\frac{\ln
q_{n+1}}{q_n}, \end{equation*}  where $\frac{p_n}{q_n}$ is the $n$-th continued fraction convergent of $\alpha.$ Recently, based on quantitative almost reducibility and Aubry duality, Avila, You and Zhou \cite{ayz} solved  the measure version of Jitomirskaya's conjecture \cite{j} on sharp phase transitions of AMO, and showed that if $\lambda>e^{\beta}$, then $H_{\lambda,\alpha,\theta}$ has AL for $a.e.$ $\theta$.  The arithmetic version was solved in \cite{JLiu}, moreover, the universal hierarchical structure of the eigenfunction was explored in \cite{JLiu}.

More recently, for any fixed Diophantine frequency, Jitomirskaya and Liu \cite{JLiu1} proved sharp transition in the arithmetics of phase between localization and singular continuous spectrum. More importantly, they uncover a new type of hierarchy for quasi-periodic operators which is called reflective-hierarchy structure. We mention that a complete characterization of eigenfunctions is also given in \cite{JLiu} which is crucial for proof of EDL in \cite{jkl}.

However, the dynamical localization results are more restricted to Diophantine frequency.
Jitomirskaya-Last  \cite{jl} proved that for $\lambda>\frac{15}{2}$, then  $H_{\lambda,\alpha,\theta}$  has semi-uniformly dynamical localization (SUDL) for $a.e.$ $\theta$:
$$\sup\limits_{t\in\R}|\langle\delta_n, e^{-itH_\theta}\delta_{\ell}\rangle| \leq \tilde{C}_{\theta,\epsilon} e^{\epsilon |\ell |}e^{-\tilde{\gamma} |n-\ell|}.$$
 Germinet and Jitomirskaya showed in \cite{gj} that for $\lambda>1$,  it has strong dynamical localization in expectation:
\begin{equation}\label{sdl}
\int_{\Omega}\sup\limits_{t}\sum\limits_{n\in\Z}(1+|n|)^q|\langle\delta_n,e^{-itH_{\theta}}\delta_0\rangle|d\theta <\infty.
\end{equation}
 Jitomirskaya and Kr\"uger \cite{jk} improved the above results by proving that AMO in fact has EDL for $\lambda>1$. However, there is no good characterization of the exponential decay rate in expectation in \cite{jk}, recently, based on the  localization method  developed in \cite{JLiu1},   Jitomirskaya, Kr\"uger and Liu \cite{jkl} proved that actually $\gamma(H)=\ln|\lambda|$.   Theorem \ref{amo}  would be seen as  an alternative proof of  their result  by different method, in fact by reducibility method as in \cite{ayz}. We  emphasize that what Jitomirskaya, Kr\"uger and Liu \cite{jkl} obtained is stronger, they are also able to prove  $$\limsup_{n\rightarrow \infty} (-\frac{1}{|n|} \ln \int_{\Omega}\sup\limits_{t\in\R}|\langle\delta_n, e^{-itH_\omega}\delta_{0}\rangle|d\mu) =\ln |\lambda|$$
based on  their fundamental work \cite{JLiu1}.  This result is surprising and has not been established for any ergodic models. We   wonder whether this is also a  general  phenomenon.

\subsection{Criterion of EDL}
We will develop a general criterion of EDL for ergodic families of operators acting on $\ell^2(\Z^d)$, the above results are actually  concrete applications of the criterion.

We first review several developments in localization theory during the past years which inspired our main ideas.
By RAGE theorem \cite{cfks}, DL implies pure point spectrum, but not vice versa, for example, one can consult the counterexample of Rio-Jitomirskaya-Last-Simon \cite{RJLS}.  Note AL means $H_{\omega}$ has a complete set of normalized eigenvectors $\{u_{\omega,m}\}_{m=1}^{\infty}$ obeying
\begin{equation}\label{eigen}
|u_{\omega,m}(n)| \leq C_{\omega,m} e^{-\gamma |n-n_{\omega,m}|},
\end{equation}
where $\gamma>0, C_{\omega,m}>0$,  and  $n_{\omega,m}\in\Z^d$ are called the centers of localization. Counterexample in  \cite{RJLS} shows that  AL is a-priori not strong enough to restrict the long time dynamics of the system. The main shortage  in \eqref{eigen} is the total freedom given to the constants $C_{\omega,m}$. Indeed, as explained in \cite{RJLS},  If  $C_{\omega,m}$ are allowed to arbitrarily grow in $m$, then, in fact, the eigenvectors can be extended over arbitrarily large length scales, possibly leading to transport arbitrarily close to the ballistic motion. Thus, to get DL, one should  specify the dependence of $C_{\omega,m}$ on $n_{\omega,m}$ and $\omega$.

Sometimes, better control for $C_{\omega,m}$ is avaliable. The breakthrough still belongs to \cite{RJLS}, they give more information on the localized eigenfunctions,  and {\it first} establish the correspondence between eigenfunctions localization and dynamical localization type result. The condition they give is called semiuniformly localized eigenvectors (SULE): for each $\epsilon>0$, there is $C_{\omega,\epsilon}$ such that
 \begin{equation}\label{eigensemi}
|u_{\omega,m}(n)| \leq C_{\omega,\epsilon} e^{\epsilon |n_{\omega,m}|}e^{-\gamma |n-n_{\omega,m}|},
\end{equation}
i.e.,  the constant $C_{\omega,m}$ has subexponential growth in the localization center $n_{\omega,m}$. They proved that SULE imply SUDL \cite{RJLS}, which means the better control of  $C_{\omega,m}$ as given in (\ref{eigensemi})  does imply DL.

In many cases, the eigenfunctions of $H_\omega$ have even better control than the one given in (\ref{eigensemi})  for most of $\omega$. To prove EDL, which is the strongest, we do require finer structure of the eigenfunctions of $H_\omega$.  However, in our criterion, we are more concentrated on the structure of the eigenfunctions. To this stage, we introduce the following concept:

\begin{Definition}\label{good}
For any $\gamma>0$, $\ell\in\Z^d$, $C(\omega)>0$, $0<C_{|\ell|}(\omega)<1$, a normalized eigenfunction\footnote{We say $u_\omega(n)$ is normalized if $\sum_n|u_\omega(n)|^2=1$.}  $u_\omega(n)$ is said to be $(\gamma,\ell,C(\omega),C_{|\ell|}(\omega))$-good, if $$|u(n)|\leq C(\omega)(e^{-\gamma |n|}+C_{|\ell|}(\omega)e^{-\gamma|n+\ell|})$$ for any $n\in\Z^d$.(Figure 1)
\end{Definition}


\begin{Remark}
Note that $C(\omega)$ is to normalize $u_\omega(n)$, $C_{|\ell|}(\omega)$ is to measure the deviation.
\end{Remark}

\begin{figure}
  \centering
  \includegraphics[width=8cm]{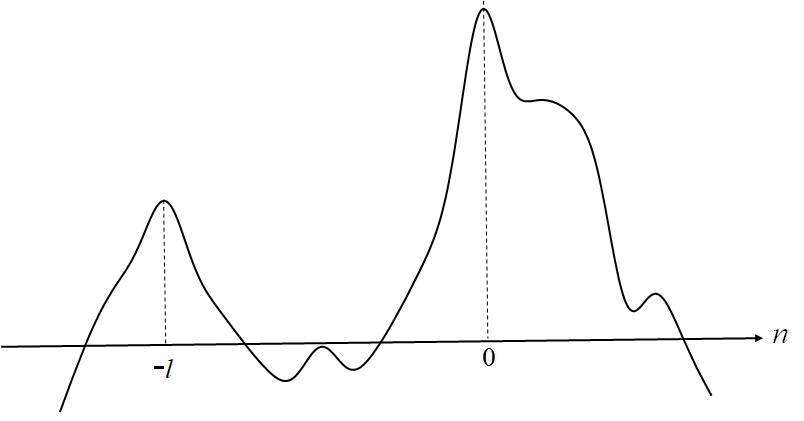} \\
  \caption{$(\gamma,\ell,C(\omega),C_{|\ell|}(\omega))$-good eigenfunction}
\end{figure}

It is obvious that every exponentially decay eigenfunction is $(\gamma,\ell,C(\omega),C_{|\ell|}(\omega))$-good for some $(\gamma,\ell,C(\omega),C_{|\ell|}(\omega))$. We remark that the constants $(\gamma,\ell,C(\omega),C_{|\ell|}(\omega))$ are not uniform in $\omega$, and  the eigenfunction is well localized if $(\ell,C_{|\ell|}(\omega))$ is small and $\gamma$ is large.  We will prove that $H_\omega$ has EDL if  each  $H_\omega$ has one sufficiently good eigenfunction for majority of $\omega$.

 Before giving the precise criterion, we first introduce a concept called {\it ergodic covariant family} which is a modification of the concept of {\it covariant family with many eigenfunctions} defined in \cite{jk2}. As proved by Bourgain and Jitomirskaya \cite{bj02},  every ergodic family of self-adjoint bounded operators with $a.e.$ pure point spectrum has a covariant spectral representation, i.e., there exist Borel measurable
functions $E(\omega) : \Omega \rightarrow \R$ and $u(\omega) :  \Omega\rightarrow \ell^2(\Z^d)$ so
that $H_{\omega}u(\omega)=E(\omega)u(\omega)$ and the functions
$\{T_{-k}u(S_k\omega),k\in\Z^d\}$  form an orthonormal basis for $\ell^2(\Z^d)$ for almost all $\omega$. It is pointed out in \cite{bj02} that covariant spectral representation have been obtained constructively for Schr\"odinger operators with Markov potentials and quasi-periodic potentials at large coupling \cite{m,Sinai}, a combination of the result in \cite{f,gjls,gk} shows that such a measurable representation always exists. Inspired by \cite{bj02} and \cite{jk2}, we have the following very natural definition.

\begin{Definition}
$(H_\omega)_{\omega\in\Omega}$  is called an ergodic covariant  family if $(H_\omega)_{\omega\in\Omega}$ is ergodic and there  exists a Borel measurable function
$$
E(\omega) : \Omega \rightarrow \R,
$$
and a normalized vector value function (not necessarily to be  measurable)
$$
u(\omega) : \Omega \rightarrow \ell^2(\Z^d),
$$
such that $H_\omega u(\omega)=E(\omega)u(\omega)$ and $\{T_{-k}u(S_k\omega), k\in\Z^d\}$ are orthonormal.
\end{Definition}

Note that the  ergodic covariant  family is weaker than covariant spectral representation, because we don't a-priori assume that the eigenfunctions are complete and $u(\omega)$ is measurable. It is also weaker than the concept of {\it covariant family with many eigenfunctions}, since we allow the operator  to have multiple eigenvalues. Moreover, it is easy to see that $(H_\omega)_{\omega\in\Omega}$ form an ergodic covariant family if $H_\omega$ has p.p. for $a.e.$ $\omega$. In fact, one can prove that the inverse is also true. i.e. if $(H_\omega)_{\omega\in\Omega}$ is an ergodic covariant family, then $H_\omega$ has p.p. for $a.e.$ $\omega$.



%
%
%

\begin{Theorem}\label{2}
Suppose that $(H_\omega)_{\omega\in\Omega}$ is an ergodic covariant family. Then $(H_{\omega})_{\omega\in\Omega}$ displays EDL  with $\gamma(H)\geq \gamma$ if the following hold\\
\begin{enumerate}
  \item{there exist $\Omega_i\nearrow\Omega$ \footnote{We say $\Omega_i\nearrow \Omega$ in the sense that $\Omega_i\subset\Omega_{i+1}$ and $\mu(\cup_i\Omega_i)=1$}, such that for any $\omega\in\Omega_i$,
 $u(\omega)$ is a $(\gamma,\ell,C_i,C_{i,m})$-good eigenfunction for some $|\ell|=m<n_i$,\\
  }
\item{
$u(\omega) :  \Omega\rightarrow \ell^2(\Z^d)$ is measurable,
}\\
\item{

$$
\sum\limits_{i=1}^\infty C^2_{i}(1+\sup_{m}C_{i,m}e^{\gamma m})\mu(\Omega^c_{i-1})<\infty.$$\\
}
\end{enumerate}

\end{Theorem}

\begin{Remark} Theorem \ref{2} is also motivated by Avila-You-Zhou \cite{ayz}, where they establish a new criterion  (Propostion 4.1) for AL.
\end{Remark}
\begin{Remark}\label{extra}
From our proof, one will see condition (1) is sufficient for proving  AL,  and the ideas of the proof of AL partially come from \cite{ayz} and \cite{jk2}.
\end{Remark}

Finally, we remark that  Jitomirskaya and her coauthors  \cite{jk, jkl} had developed a method to prove EDL, and the key of their approach is to exploit the orthogonality of the $u_{\omega,m}$ when estimating $\sum_{n_{\omega,m}=n}|u_{\omega,m}(\ell)|^2$. As reader can see, their method is completely different from us. 
We emphasize that the criterion in Theorem \ref{2} applies to general  ergodic  family of operators. In this paper, we only give  its applications  to quasi-periodic Schr\"odinger operators.

\subsection{Plan of the proof}  Our approach is from the perspective of dynamical systems, and  is based on quantitative almost reducibility (c.f. section \ref{aglobal}). The philosophy is that nice quantitative almost reducibility would imply nice spectral applications.
This approach, which started from the  pioneer work  of Eliasson \cite{E92},  has been proved to be very fruitful \cite{aj1,AK06, ayz, ayz1,LYZZ}.

In our case, we will first prove Theorem \ref{2} which can be viewed as a generalization of the criterion for Anderson localization in \cite{ayz}.
Second,  we apply the criterion to quasi-periodic  models. For this purpose, we need  a good  control of $C(\omega), C_{|\ell|}(\omega)$ and $l$ for $u(\omega)$. The control will be given by  Aubry duality and the quantitative almost reducibility of the dual cocycles.
More precisely, we need to analyze the behavior of the dual Schr\"odinger cocycles with Diophantine rotation numbers. We know, in this case,   the cocycles are reducible to a constant cocycle with elliptic eigenvalues. The key point is to prove the exponentially decay of the off-diagonal element of the constant matrix and  fine structure of the transformations. Furthermore, to get the optimal decay rate of the eigenfunctions, we need a strong almost reducibility result, i.e. the cocycle is almost reducible in a fixed band, moreover, we need the band arbitrary close to the initial band.

We mention that the KAM method only works for the local regime. In order to pass it to the global regime, we need Avila's global theory of analytic $SL(2,\R)$-cocycles \cite{avila0}. Moreover, in order to have a uniform  control on the conjugacies with respect to subcritical energy, we shall perform some compactness argument, the key still follows from Avila's global theory: openness of the subcritical cocycle and openness of almost reducibility.

\section{Preliminaries}

Recall $sl(2,\R)$ is the set of $2\times 2$ matrices  of the form
$$\left(
\begin{array}{ccc}
 x &  y+z\\
 y-z &  -x
 \end{array}\right)$$
 where $x,y,z\in \R.$ It is
 isomorphic to $su(1,1)$, the group of matrices of the form
$$\left(
\begin{array}{ccc}
 i t &  \nu\\
\bar{ \nu} &  -i t
 \end{array}\right)$$
 with $t\in \R$, $\nu\in \C$. The isomorphism between $sl(2,\R)$ and
 $su(1,1)$ is given by $B\rightarrow M B M^{-1}$ where
$$M=\frac{1}{2i}\left(
\begin{array}{ccc}
 1 &  -i\\
 1 &  i
 \end{array}\right).$$
 Direct calculation shows that
 $$M\left(
\begin{array}{ccc}
 x &  y+z\\
 y-z &  -x
 \end{array}\right)M^{-1}=\left(
\begin{array}{ccc}
 i z &  x-iy\\
x+iy &  -i z
 \end{array}\right).$$

\subsection{Continued Fraction Expansion} Let $\alpha\in (0,1)\backslash\Q$, $a_0:=0$ and $\alpha_0:=\alpha$. Inductively, for $k\geq 1$, we define
$$
a_k:=[\alpha_{k-1}^{-1}], \  \ \alpha_k=\alpha_{k_1}^{-1}-a_k.
$$
Let $p_0:=0$, $p_1:=1$, $q_0:=1$, $q_1:=a_1$. We define inductively $p_k:=a_kp_{k-1}+p_{k-2}$, $q_k:=a_kq_{k-1}+q_{k-2}$. Then $q_n$ are the  denominators of the best rational approximations of $\alpha$ since we have $\|k\alpha\|\geq \|q_{n-1}\alpha\|$ for all $k$ satisfying $\forall 1\leq k< q_n$, and
$$
\frac{1}{2q_{n+1}}\leq \|q_n\alpha\|_{\R/\Z}\leq \frac{1}{q_{n+1}}.
$$

\subsection{Quasi-periodic cocycles}

Given $A \in C^\omega(\T^d,{\rm SL}(2,\C))$ and rationally independent $\alpha\in\R^d$, we define the quasi-periodic \textit{cocycle} $(\alpha,A)$:
$$
(\alpha,A)\colon \left\{
\begin{array}{rcl}
\T^d \times \C^2 &\to& \T^d \times \C^2\\[1mm]
(x,v) &\mapsto& (x+\alpha,A(x)\cdot v)
\end{array}
\right.  .
$$
The iterates of $(\alpha,A)$ are of the form $(\alpha,A)^n=(n\alpha,  \mathcal{A}_n)$, where
$$
\mathcal{A}_n(x):=
\left\{\begin{array}{l l}
A(x+(n-1)\alpha) \cdots A(x+\alpha) A(x),  & n\geq 0\\[1mm]
A^{-1}(x+n\alpha) A^{-1}(x+(n+1)\alpha) \cdots A^{-1}(x-\alpha), & n <0
\end{array}\right.    .
$$
The {\it Lyapunov exponent} is defined by
$\displaystyle
L(\alpha,A):=\lim\limits_{n\to \infty} \frac{1}{n} \int_{\T^d} \ln \|\mathcal{A}_n(x)\| dx
$.

The cocycle $(\alpha,A)$ is {\it uniformly hyperbolic} if, for every $x \in \T^d$, there exists a continuous splitting $\C^2=E^s(x)\oplus E^u(x)$ such that for every $n \geq 0$,
$$
\begin{array}{rl}
|A_n(x) \, v| \leq C e^{-cn}|v|, &  v \in E^s(x),\\[1mm]
|A_n(x)^{-1}   v| \leq C e^{-cn}|v|, &  v \in E^u(x+n\alpha),
\end{array}
$$
for some constants $C,c>0$.
This splitting is invariant by the dynamics, i.e.,
$$A(x) E^{*}(x)=E^{*}(x+\alpha), \quad *=``s"\;\ {\rm or} \;\ ``u", \quad \forall \  x \in \T^d.$$

Assume that $A \in C (\T^d, {\rm SL}(2, \R))$ is homotopic to the identity. It induces the projective skew-product $F_A\colon \T^d \times \mathbb{S}^1 \to \T^d \times \mathbb{S}^1$ with
$$
F_A(x,w):=\left(x+\a,\, \frac{A(x) \cdot w}{|A(x) \cdot w|}\right),
$$
which is also homotopic to the identity.
Thus we can lift $F_A$ to a map $\tF_A\colon \T^d \times \R \to \T^d \times \R$ of the form $\tF_A(x,y)=(x+\alpha,y+\psi_x(y))$, where for every $x \in \T^d$, $\psi_x$ is $\Z$-periodic.
The map $\psi\colon\T^d \times \T  \to \R$ is called a {\it lift} of $A$. Let $\mu$ be any probability measure on $\T^d \times \R$ which is invariant by $\widetilde{F}_A$, and whose projection on the first coordinate is given by Lebesgue measure.
The number
$$
\rho(\alpha,A):=\int_{\T^d \times \R} \psi_x(y)\ d\mu(x,y) \ {\rm mod} \ \Z
$$
 depends  neither on the lift $\psi$ nor on the measure $\mu$, and is called the \textit{fibered rotation number} of $(\alpha,A)$ (see \cite{H,JM} for more details).

Given $\theta\in\T^d$, let $
R_\theta:=
\begin{pmatrix}
\cos2 \pi\theta & -\sin2\pi\theta\\
\sin2\pi\theta & \cos2\pi\theta
\end{pmatrix}$.
If $A\colon \T^d\to{\rm PSL}(2,\R)$ is homotopic to $\theta \mapsto R_{\frac{\la n, \theta\ra}{2}}$ for some $n\in\Z^d$,
then we call $n$ the {\it degree} of $A$ and denote it by $\deg A$.
The fibered rotation number is invariant under real conjugacies which are homotopic to the identity. More generally, if $(\alpha,A_1)$ is conjugated to $(\alpha, A_2)$, i.e., $B(\cdot+\alpha)^{-1}A_1(\cdot)B(\cdot)=A_2(\cdot)$, for some $B \colon \T^d\to{\rm PSL}(2,\R)$ with ${\rm deg} B=n$, then
\begin{equation}\label{rotation number}
\rho(\alpha, A_1)= \rho(\alpha, A_2)+ \frac{\la n,\alpha \ra}2.
\end{equation}

%

A typical  example is given by the so-called \textit{Schr\"{o}dinger cocycles} $(\alpha,S_E^{V})$, with
$$
S_E^{V}(\cdot):=
\begin{pmatrix}
E-V(\cdot) & -1\\
1 & 0
\end{pmatrix},   \quad E\in\R.
$$
Those cocycles were introduced because it is equivalent to the eigenvalue equation $H_{V, \alpha, \theta}u=E u$. Indeed, any formal solution $u=(u_n)_{n \in \Z}$ of $H_{V, \alpha, \theta}u=E u$ satisfies
$$
\begin{pmatrix}
u_{n+1}\\
u_n
\end{pmatrix}
= S_E^V(\theta+n\alpha) \begin{pmatrix}
u_{n}\\
u_{n-1}
\end{pmatrix},\quad \forall \  n \in \Z.
$$
The spectral properties of $H_{V,\alpha,\theta}$ and the dynamics of $(\alpha,S_E^V)$ are closely related by the well-known fact:
 $E\in \Sigma_{V,\alpha}$ if and only if $(\alpha,S_E^{V})$ is \textit{not} uniformly hyperbolic. Throughout the paper, we will denote $L(E)=L(\alpha,S_E^{V})$  and $\rho(E)=\rho(\alpha,S_E^{V})$ for short.

The cocycles $(\alpha,A)$ is $C^\omega$-reducible, if it can be $C^\omega$-conjugated to a constant cocycle.

%
%

\subsection{Global theory of one-frequency Schr\"odinger operators} \label{aglobal}Let us make a short review of Avila's global theory of one frequency $SL(2,\R)$-cocycles \cite{avila0}. Suppose that $A\in C^\omega(\T, SL(2,\R))$ admits a holomorphic extension to $\{|\Im z|<h\}$. Then for $|\epsilon|<h$, we define $A_\epsilon\in C^\omega(\T,SL(2,\C))$ by $A_\epsilon(\cdot)=A(\cdot+i\epsilon)$. The cocycles which are not uniformly hyperbolic are classified into three classed: subcritical, critical, and supercritical. In particular, $(\alpha,A)$ is said to be subcritical if there exists $h>0$ such that $L(\alpha,A_\epsilon)=0$ for $|\epsilon|<h$.

A cornerstone in Avila's global theory is the ``Almost Reducibility Conjecture" (ARC), which says that $(\alpha,A)$ is almost reducible if it is subcritical. Recall that the cocycle $(\alpha,A)$ is said to be almost reducible if there exists $h_*>0$, and a sequence $B_n\in C_{h_*}^\omega(\T,PSL(2,R))$ such that $B^{-1}_n(\theta+\alpha)A(\theta)B(\theta)$ converges to constant uniformly in $|\Im\theta|<h_*$. The complete solution of ARC was recently given by Avila \cite{avila1,avila2}, in the case $\beta(\alpha)=0$, it is the following:
\begin{Theorem} \label{arc-conjecture}
For $\alpha\in \R\backslash\Q$ with $\beta(\alpha)=0$, and $A\in C^\omega(\T,SL(2,\R))$, $(\alpha,A)$ is almost reducible if it is subcritical.
\end{Theorem}

\subsection{$\Z^d$-action and ergodicity}
Generally speaking, an action of a group is a formal way of interpreting the manner in which the elements of the group correspond to transformations of some space, in a way that preserves the structure of that space.  An important and special group is $\Z^d$ additive group.

Specially, the group $S=\{S_n,n\in\Z^d\}$, a family of measure-preserving one to one transformations of $\Omega$ such that $S_{m+n}=S_mS_n$,  can be seen as a $\Z^d$-action on $\Omega$. The measurable dynamical system $(\Omega,d\mu,S)$ is called ergodic if it holds that  for any $f\in L^1(\Omega,d\mu)$ and for $\mu$ $a.e.$ $\omega\in \Omega$,
\begin{equation}\label{ergodic}
\lim\limits_{N\rightarrow \infty}\frac{1}{(2N)^d}\sum\limits_{|n|\leq N}f(S_n\omega)=\int_{\Omega}f d\mu.
\end{equation}
For more details, one can consult \cite{ps}.

\subsection{Aubry duality}
Suppose that the quasi-periodic Schr\"{o}dinger operator
\begin{equation}
(H_{\lambda^{-1} V,\alpha,\theta}u)_n=u_{n+1}+u_{n-1}+\lambda^{-1} V(\theta+n\alpha)x_{n}, n\in \Z.
\end{equation}
has an analytic quasi-periodic Bloch wave $u_n=e^{2\pi i n\varphi}\overline{\psi}(\phi+n\alpha)$ for some $\overline{\psi}\in C^\omega(\T^d,\C)$ and $\varphi \in {[}0,1{)}$. It is easy to see that the Fourier coefficients of $\overline{\psi}(\theta)$ is an eigenfunction of  the following long range operator:
\begin{equation}
(L_{V, \lambda, \alpha,\theta}u)_m=\sum_{k\in \Z^d}V_k u_{m-k}+2\lambda\cos 2\pi(\varphi+\langle m,\alpha\rangle)u_m, m\in \Z^d.
\end{equation}
$L_{V, \lambda, \alpha,\theta}$ is called the dual operator of $H_{\lambda^{-1} V,\alpha,\theta}$.

\section{Criterion of EDL}
In this section, we give the proof of Theorem \ref{2}. Denote $E_k(\omega)=E(S_k\omega)$, $u_k(\omega)=T_{-k}u(S_k\omega)$ and denote $u_k(n,\omega)$ the $n$-th component of $u_k(\omega)$. Then we have
\begin{Lemma}
If $(H_\omega)_{\omega\in\Omega}$ is an ergodic family of operators, then $E_k(\omega)$ is an eigenvalue of $H_\omega$ with eigenfunction $u_k(\omega)$.
\end{Lemma}
\begin{pf}
we have
$$ H_\omega T_{-k}u(S_{k}\omega)=T_{-k}H_{S_k\omega} u(S_k\omega)=T_{-k}E(S_k\omega)u(S_k\omega)=E(S_k\omega)T_{-k}u(S_{k}\omega),$$
 which means that $E_k(\omega)$ is an eigenvalue of $H_\omega$ with eigenfunction $u_k(\omega)$.
\end{pf}

Now by ergodicity, there exists $\tilde{\Omega}$ with $\mu(\tilde{\Omega})=1$, such that for any $\omega\in \tilde{\Omega}$ and the measurable set $\Omega_i$ defined in Theorem \ref{2},
\begin{equation}\label{condition}
\#\{k|S_k\omega\in\Omega_i,|k|\leq \mu(\Omega_i)N\}\geq \mu^3(\Omega_i)(2N)^d
\end{equation} holds for $N>N_0$ where $N_0$ depends on $\omega$ and $\Omega_i$.
We denote
\begin{flalign*}
K_{\gamma,N,i}(\omega)=\{k||k|\leq \mu(\Omega_i)N, S_k\omega\in \Omega_i\},
\end{flalign*}

The following lemma is a quantitative description of $u_{k}(\omega)$ for $k\in K_{\gamma,N,i}(\omega)$.
\begin{Lemma}\label{6}
For any $k\in K_{\gamma,N,i}(\omega)$ and $p\in\Z^d$,
we have
\begin{equation*}
\sum\limits_{|n-p|\leq N}|u_{k}(n,\omega)|^2\geq 1-e^{-\frac{\gamma(1-\mu(\Omega_i))}{2} N},
\end{equation*}
provided that
$$
N>N_0':=\max\limits_{m\leq n_i}\{\frac{2n_i+2|p|+2\ln C_i(1+C_{i,m})}{\gamma(1-\mu(\Omega_i))},(d+1)!(\frac{2}{\gamma(1-\mu(\Omega_i))})^{d+1}\}.
$$
\end{Lemma}
\begin{pf}
Note that for any $k\in K_{\gamma,N,i}(\omega)$, we have $S_k\omega\in \Omega_i$, by condition (1), $u(S_k\omega)$ is a $(\gamma,\ell,C_i,C_{i,m})$-good eigenfunction for some $|\ell|=m<n_i$, since $u_k(\omega)=T_{-k}u(S_k\omega)$, we have
$$
|u_k(n,\omega)|\leq C_i(e^{-\gamma|n+k|}+C_{i,m}e^{-\gamma|n+k+\ell|}),
$$
thus
\begin{align*}
\sum\limits_{|n-p|\geq N}|u_{k}(n,\omega)|^2&\leq \sum\limits_{|n-p|\geq N}C_i(e^{\gamma(|k|+|p|)}e^{-\gamma|n-p|}+C_{i,m}e^{\gamma(|k|+|n_i|+|p|)}e^{-\gamma|n-p|})\\
&\leq \sum\limits_{j\geq N}j^dC_i(1+C_{i,m})e^{\gamma(|k|+|n_i|+|p|)}e^{-\gamma j}\\
&\leq \sum\limits_{j\geq N}j^dC_i(1+C_{i,m})e^{\gamma(|n_i|+|p|)}e^{-\gamma(1-\mu(\Omega_i)) j}.
\end{align*}
By the definition of $N_0'$, direct computation shows for $N>N_0'$, we have
\begin{equation*}
\sum\limits_{|n-p|\leq N}|u_{k}(n,\omega)|^2\geq 1-e^{-\frac{\gamma(1-\mu(\Omega_i))}{2} N}>0,
\end{equation*}

\end{pf}

Once we have this,  we denote
\begin{align*}
\mathcal{E}_{\gamma,N,i}(\omega)=\{E_k(\omega)||k|\leq \mu(\Omega_i)N, S_k\omega\in \Omega_i\}.
\end{align*}
$$
\mathcal{E}_{\gamma}(\omega)=\bigcup\limits_{i\in\N}\bigcup\limits_{N\in\N}\mathcal{E}_{\gamma,N,i}(\omega).
$$
We remark that $E_k(\omega)$ might be multiple eigenvalues, since we didn't assume $E_m(\omega)\neq E_n(\omega)$ for $m\neq n$. However, by our assumption, all the eigenfunctions of $E(S_k\omega)$ are orthonomal, which is quite useful for our following computations.

Let $P^\omega_{\Delta}$ be the spectral projection of $H_\omega$ onto the set $\Delta$, and for any $\delta_n\in\ell^2(\Z^d)$, it induces the spectral measure:
$$
\langle P_\Delta^\omega\delta_n,\delta_n\rangle=\mu_{\delta_n,\omega}(\Delta).
$$
For $N>\max\{N_0,N_0'\}$, by the spectral theorem, \eqref{condition} and Lemma \ref{6}, we have
\begin{align*}
\frac{1}{(2N)^d}\sum\limits_{|n-p|\leq N}|\mu_{\delta_n,\omega}(\mathcal{E}_{\gamma}(\omega))|&>\frac{1}{(2N)^d}\sum\limits_{|n-p|\leq N}<P^\omega_{\mathcal{E}_{\gamma,N,i}(\omega)}\delta_n,\delta_n>\\
&>\frac{1}{(2N)^d}\sum\limits_{|n-p|\leq N}\sum\limits_{k\in K_{\gamma,N,i}(\omega)}|u_k(n,\omega)|^2\\
&>\frac{1}{(2N)^d}\#K_{\gamma,N,i}(\omega)(1-e^{-\frac{\gamma(1-\mu(\Omega_i))}{2} N})\\
&\geq\mu(\Omega_i)^{3}(1-e^{-\frac{\gamma(1-\mu(\Omega_i))}{2} N}).
\end{align*}
Note the second inequality uses the fact that all the eigenfunctions of $E(S_k\omega)$ are orthonormal. Since $\mathcal{E}_{\gamma}(\omega)=\mathcal{E}_{\gamma}(S_n\omega)$ where $S_n$ is the measure-preserving transformation defined in \eqref{ergodic}, we can rewrite the above inequalities as
$$
\frac{1}{(2N)^d}\sum\limits_{|n|\leq N}|\mu_{\delta_p,S_n\omega}(\mathcal{E}(S_n\omega))|>\mu(\Omega_i)^{3}(1-e^{-\frac{\gamma(1-\mu(\Omega_i))}{2} N}).
$$

Let $N$ goes to $\infty$, thus by ergodicity. we have
$$
\int_\Omega|\mu_{\delta_p,\omega}(\mathcal{E}_{\gamma}(\omega))|d\omega\geq \mu(\Omega_i)^{3},
$$
Let $i\rightarrow \infty$, it follows that
\begin{equation}\label{simple}
\int_\Omega|\mu_{\delta_p,\omega}(\mathcal{E}_{\gamma}(\omega))|d\omega=1.
\end{equation}

Thus for $\mu$-$a.e.$ $\omega\in \Omega$, $\mu_{\delta_p,\omega}=\mu_{\delta_p,\omega}^{pp}$ with support $\mathcal{E}_{\gamma}(\omega)$, by passing to the countable intersection of these set of $\omega$, we ultimately obtain a full measure set for which $H_\omega$ has pure point spectrum with eigenfunctions $(u_k(\omega))_{k\in\Z^d}$. Since  $(u_k(\omega))_{k\in\Z^d}$ all decays exponentially, $H_\omega$ has AL for $a.e.$ $\omega$.

Next we prove EDL of the ergodic covariant family $(H_\omega)_{\omega\in\Omega}$. For any $i,m\in \N$, we inductively define
$$
\Omega_i^0=\{\mbox{$\omega\in\Omega_i|u(\omega)$ is a $(\gamma,0,C_i,C_{i,0})$-good eigenfunction}\},
$$
\begin{align*}
\Omega_i^m=&\{\mbox{$\omega\in\Omega_i|\omega\notin \bigcup\limits_{j=1}^{m-1}\Omega_i^j,u(\omega)$ is a $(\gamma,\ell,C_i,C_{i,m})$}\\
 & \mbox{-good eigenfunction for some $|\ell|=m$}\},
\end{align*}
by condition (1), $\Omega_i=\bigcup\limits_{m=0}^{n_i}\Omega_i^m$. Since $u(\omega)$ is measurable, we have $\{\Omega_i^m\}_{m=0}^{n_i}$ are all measurable sets, and it is obvious that $\Omega_i^{m_1}\cap\Omega_i^{m_2}=\emptyset$ for any $m_1\neq m_2$.

Now we need the following  standard computation:
\begin{equation*}
\int_\Omega\langle\delta_p,e^{-itH(\omega)}\delta_q\rangle d\mu(\omega)=\sum\limits_{k\in\Z^d}\int_\Omega e^{-itE_k(\omega)}\overline{u_k(p,\omega)}u_k(q,\omega)d\mu(\omega),
\end{equation*}
thus
\begin{equation*}
\int_\Omega|\langle\delta_p,e^{-itH(\omega)}\delta_q\rangle| d\mu(\omega)\leq\sum\limits_{k\in\Z^d}\int_\Omega |\overline{u_k(p,\omega)}u_k(q,\omega)|d\mu(\omega).
\end{equation*}
Since $\{S_k\}_{k\in\Z^d}$ is a family of measure preserving transformation, we have
\begin{align*}
&\sum\limits_{k\in\Z^d}\int_\Omega |\overline{u_k(p,\omega)}u_k(q,\omega)|d\mu(\omega)\\ =&\sum\limits_{k\in\Z^d}\int_\Omega |\overline{u_k(p,S_{-k}\omega)}u_k(q,S_{-k}\omega)|d\mu(\omega)\\
=&\sum\limits_{k\in\Z^d}\int_\Omega |\overline{u_0(p-k,\omega)}u_0(q-k,\omega)|d\mu(\omega)\\
=&\sum\limits_{k\in\Z^d}\sum\limits_{i=1}^\infty\int_{\Omega_i\backslash\Omega_{i-1}}|\overline{u_0(p-k,\omega)}u_0(q-k,\omega)|d\mu(\omega)\\
=&\sum\limits_{i=1}^{i_0} \sum\limits_{k\in\Z^d}\int_{\bigcup\limits_{i=0}^{n_i}\Omega_i^m\backslash\Omega_{i-1}}|\overline{u_0(p-k,\omega)}u_0(q-k,\omega)|d\mu(\omega) \\ +&\sum\limits_{k\in\Z^d}\int_{\Omega_{i_0}^c}|\overline{u_0(p-k,\omega)}u_0(q-k,\omega)|d\mu(\omega).
\end{align*}
The second quality holds because of $u_k(n,S_{-k}\omega)=T_ku(n,\omega)=u(n-k,\omega)$.

Assume $i_0$ is the smallest integer $i$ such that $2\mu(\Omega^c_i)<e^{-\gamma|p-q|}$, by Bessel's inequality, we have
$$
\sum\limits_{k\in\Z^d}\int_{\Omega_{i_0}^c}|\overline{u_0(p-k,\omega)}u_0(q-k,\omega)|d\mu(\omega)\leq \mu(\Omega_{i_0}^c) \leq e^{-\gamma|p-q|}.
$$
On the other hand, by the definition of $\Omega_i^m$ and direct computation, one has
\begin{align*}
&\ \ \ \ \int_{\Omega_i^m\backslash\Omega_{i-1}}|\overline{u_0(p-k,\omega)}u_0(q-k,\omega)| d\mu(\omega)\\ \nonumber
&<C_i^2\mu(\Omega_i^m\backslash\Omega_{i-1})\big(e^{-\gamma(|p-k|+|q-k|)}+C_{i,m}^2e^{-\gamma(|p-k+\ell|+|q-k+\ell|)}\\ \nonumber
&\ \ +C_{i,m}(e^{-\gamma(|p-k|+|q-k+\ell|)}+e^{-\gamma(|p-k+\ell|+|q-k|)})\big),\nonumber
\end{align*}
where $|\ell|=m$.

To complete the proof, we need the following simple observations:
\begin{Lemma}\label{criterionlemma2}We have the following inequality:
\begin{eqnarray}
 \label{crelemma-1}\sum\limits_{k\in\Z^d}e^{-\gamma(|p-k|+|q-k|)} &\leq&  (2C(\gamma)+d+|p-q|)^de^{-\gamma|p-q|},\\
 \label{crelemma-2}\sum\limits_{k\in\Z^d}e^{-\gamma(|p-k+\ell|+|q-k+\ell|)} &\leq&  (2C(\gamma)+d+|p-q|)^de^{-\gamma|p-q|},
 \end{eqnarray}
\begin{eqnarray}
\nonumber &\sum\limits_{k\in\Z^d}(e^{-\gamma(|p-k|+|q-k+\ell|)}+e^{-\gamma(|p-k+\ell|+|q-k|)}) \\  \label{crelemma-3}  &\leq  2e^{\gamma |\ell|}(2C(\gamma)+d+|p-q|)^de^{-\gamma|p-q|}.
\end{eqnarray}

\begin{pf} Let $p=(p_1,\cdots,p_d)$, $q=(q_1,\cdots,q_d)$, $k=(k_1,\cdots,k_d)$, $\ell=(\ell_1,\cdots,\ell_d)$, and assume that $p_1\leq q_1$, let $\Lambda_1=\{k_1\in\Z|k_1<p_1\}$, $\Lambda_2=\{k_1\in\Z|k_1>q_1\}$, it is obvious that
$$
\sum\limits_{k_1\in\Lambda_1}e^{-\gamma(|p_1-k_1|+|q_1-k_1|)}\leq \sum\limits_{k_1=-\infty}^{p_1}e^{-\gamma|q_1-k_1|}\leq C(\gamma)e^{-\gamma|p_1-q_1|},
$$
$$
\sum\limits_{k_1\in\Lambda_2}e^{-\gamma(|p_1-k_1|+|q_1-k_1|)}\leq \sum\limits_{k_1=q_1}^{\infty}e^{-\gamma|p_1-k_1|}\leq C(\gamma)e^{-\gamma|p_1-q_1|},
$$
since $\#(\Lambda_1\cup\Lambda_2)^c= |p_1-q_1+1|$ and for $k_1\in(\Lambda_1\cup\Lambda_2)^c$, $$|p_1-k_1|+|q_1-k_1|= |p_1-q_1|,$$ we have
$$
\sum\limits_{k_1\in(\Lambda_1\cup\Lambda_2)^c}e^{-\gamma(|p_1-k_1|+|q_1-k_1|)}= \sum\limits_{k_1\in(\Lambda_1\cup\Lambda_2)^c}e^{-\gamma|p_1-q_1|}= |p_1-q_1+1|e^{-\gamma|p_1-q_1|},
$$
therefore
$$
\sum\limits_{k_1\in\Z}e^{-\gamma(|p_1-k_1|+|q_1-k_1|)}\leq (2C(\gamma)+|p_1-q_1+1|)e^{-\gamma|p_1-q_1|}.
$$

By the definition of $l^1$ norm, we have $$
|p-k|+|q-k|=\sum\limits_{i=1}^d|p_i-k_i|+|q_i-k_i|,
$$
therefore
\begin{align*}
\sum\limits_{k\in\Z^d}e^{-\gamma(|p-k|+|q-k|)}&\leq e^{-\gamma|p-q|}\prod\limits_{i=1}^d(2C(\gamma)+|p_i-q_i+1|)\nonumber \\
&\leq (2C(\gamma)+d+|p-q|)^d e^{-\gamma|p-q|},
\end{align*}
this finishes the proof of $(\ref{crelemma-1})$. The proof of $(\ref{crelemma-2})$ is similar, we omit the details.
On the other hand, we have $e^{-\gamma|p-k+\ell|}\leq e^{\gamma|\ell|}e^{-\gamma|p-k|}$, therefore by $(\ref{crelemma-1})$, we immediately get the proof of $(\ref{crelemma-3})$. \end{pf}
\end{Lemma}
By Lemma \ref{criterionlemma2}, we have
\begin{align*}
&\ \ \ \ \sum\limits_{i=1}^{i_0} \sum\limits_{k\in\Z^d}\int_{\bigcup\limits_{i=0}^{n_i}\Omega_i^m\backslash\Omega_{i-1}}|\overline{u_0(p-k,\omega)}u_0(q-k,\omega)|d\mu(\omega)\\
&\leq\sum\limits_{i=0}^{i_0}\sum\limits_{m=0}^{n_i}C_i^2(1+C_{i,m}^2+2C_{i,m}e^{\gamma m})\mu(\Omega_i^{m}\backslash\Omega_{i-1})(2C(\gamma)+d+|p-q|^d)e^{-\gamma|p-q|}\\
&\leq 2(2C(\gamma)+d+|p-q|^d) \sum\limits_{i=0}^{i_0}C_i^2(1+\sup_{m}C_{i,m}e^{\gamma m})\mu(\Omega^c_{i-1})e^{-\gamma|p-q|},
\end{align*}
therefore by condition (3), it implies that
\begin{align*}
|\int_\Omega\langle\delta_p,e^{-itH(\omega)}\delta_q\rangle d\mu(\omega)|&\leq C(2C(\gamma)+d+|p-q|^d)e^{-\gamma|p-q|}.
\end{align*}
By the definition of $\gamma(H)$, one thus has $\gamma(H)\geq \gamma$.

\qed

\section{Good eigenfunctions via reducibility}
In this section, we give a method for obtaining $(\gamma, \ell,C,C_{|\ell|})$-good eigenfunctions from the  quantitative reducibility of the dual Schr\"odinger cocycles. The key point is that we not only need quantitative estimates of the transformation, but also the structure of it.

\begin{Proposition}\label{Upper bound}
Let $\alpha\in \T^d$, $0<h'<h$, $V\in C_h^{\omega}(\T^d,\R).$ Assume that $\rho(E)$ is not rational w.r.t. $\alpha$ and
$(\alpha,S_E^{\lambda^{-1} V})$ is reducible, i.e.  there exist
 $B\in C^\omega_{h'}(\T^d,PSL(2,\R))$ and $A\in SL(2,\R)$ such that
$$
B^{-1}(\theta+\alpha)S_E^{\lambda^{-1} V}(\theta)B(\theta)=A=M^{-1}exp \left(
\begin{array}{ccc}
 i t &  \nu\\
\bar{ \nu} &  -i t
 \end{array}\right)M.
$$
Moreover, assume that the conjugacy $B(\theta)$ can be written as
$$
B(\theta)=\tilde{B}(\theta)R_{\frac{\langle \ell,\theta\rangle}{2}}e^{Y(\theta)},
$$
with $\|Y\|_{h'} \leq \frac{1}{2}.$
Then the long-range operator $L_{V,\lambda,\alpha,\rho(E)}$ has a $(2\pi h',\pm\ell,C,C_{|\ell|})$-good eigenfunction with \begin{eqnarray*}
C &=& 8\|\tilde{B}\|_{h'}^{4}\\
C_{|\ell|}&=& \min\{1,\|Y\|_{h'}+\frac{2|\nu|}{\|2\rho(E)-\langle \ell,\alpha\rangle-\langle \ell_0,\alpha\rangle\|_{\R/\Z}}\},
\end{eqnarray*}
where $\ell_0=\deg{\tilde{B}}$.
\end{Proposition}

\begin{Remark}
In Proposition \ref{Upper bound}, Diophantine condition on $\alpha$ is not assumed.
\end{Remark}

\begin{pf}
By the assumption,  the conjugacy $e^{Y(\theta)}$ is close to identity, thus it has zero degree, and then $\deg{B}= \ell+ \deg{\tilde{B}}=\ell +\ell_0 $. Now since $(\alpha,S_E^{\lambda^{-1} V})$ is reducible  to $(\alpha,A)$,  by \eqref{rotation number}, we have
$$\rho(\alpha,A)=\rho(E)-\frac{\langle \ell+\ell_0,\alpha\rangle}{2}.$$
Since $\rho(E)$ is not rational w.r.t. $\alpha$, we furthermore have
 $$spec{A}=\{e^{2\pi i(\rho(E)-\frac{\langle \ell+\ell_0,\alpha\rangle}{2})},e^{-2\pi i(\rho(E)-\frac{\langle \ell+\ell_0,\alpha\rangle}{2})}\}.$$

Since $A\in SL(2,\R)$, there exists  unitary $U\in SL(2,\C)$ such that
$$
UMAM^{-1}U^{-1}=\begin{pmatrix}e^{2\pi i(\rho(E)-\frac{\langle \ell+\ell_0,\alpha\rangle}{2})}&c\\0&e^{2\pi i(-\rho(E)+\frac{\langle \ell+\ell_0,\alpha\rangle}{2})}\end{pmatrix}.
$$
Let $B_1(\theta)=\tilde{B}(\theta)R_{\frac{\langle \ell,\theta\rangle}{2}}e^{Y(\theta)}M^{-1}U^{-1}=\begin{pmatrix}b_{11}(\theta)&b_{12}(\theta)\\b_{21}(\theta)&b_{22}(\theta)\end{pmatrix}$, we have
$$
B_1^{-1}(\theta+\alpha)S_E^{\lambda^{-1} V}(\theta)B_1(\theta)=\begin{pmatrix}e^{2\pi i(\rho(E)-\frac{\langle \ell+\ell_0,\alpha\rangle}{2})}&c\\0&e^{2\pi i(-\rho(E)+\frac{\langle \ell+\ell_0,\alpha\rangle}{2})}\end{pmatrix}.
$$
This implies that
\begin{equation}\label{f2}
b_{11}(\theta)=e^{2\pi i(\rho(E)-\frac{\langle \ell+\ell_0,\alpha\rangle}{2})}b_{21}(\theta+\alpha),
\end{equation}
and consequently,
\begin{align*}
&(E-\lambda^{-1} V(\theta))b_{11}(\theta)\\ \nonumber
=&b_{11}(\theta-\alpha)e^{-2\pi i(\rho(E)-\frac{\langle \ell+\ell_0,\alpha\rangle}{2})}+b_{11}(\theta+\alpha)e^{2\pi i(\rho(E)-\frac{\langle \ell+\ell_0,\alpha\rangle}{2})}.
\end{align*}


We denote $z_{11}(\theta)=e^{-2\pi i\frac{\langle \ell+\ell_0,\theta\rangle}{2}}b_{11}(\theta)$, then one has
\begin{align}\label{4.1}
(E-\lambda^{-1} V(\theta))z_{11}(\theta)=z_{11}(\theta-\alpha)e^{-2\pi i\rho(E)}+z_{11}(\theta+\alpha)e^{2\pi i\rho(E)}.
\end{align}
Taking the Fourier transformation of \eqref{4.1}, we have
$$
\sum\limits_{k\in\Z^d}\hat{z}_{11}(n-k)V_k+2\lambda\cos2\pi(\rho(E)+\langle n,\alpha\rangle)\hat{z}_{11}(n)=\lambda E\hat{z}_{11}(n),
$$
i.e.  $\{\hat{z}_{11}(n),n\in\Z^d\}$ is an eigenfunction of the long-range operator $L_{V,\lambda,\alpha,\rho(E)}$.

Now we give the estimate of the eigenfunction $\{\hat{z}_{11}(n),n\in\Z^d\}$. For convenience, we denote
 \begin{equation*}  e^{MY(\theta)M^{-1}}=Me^{Y(\theta)}M^{-1}=\begin{pmatrix}y_{11}(\theta)&y_{12}(\theta)\\y_{21}(\theta)&y_{22}(\theta)\end{pmatrix}, \end{equation*}
 $$U^{-1}=\begin{pmatrix}u_{11}&u_{12}\\u_{21}&u_{22}\end{pmatrix},\quad \tilde{B}(\theta)M^{-1}=\begin{pmatrix}\tilde{b}_{11}(\theta)&\tilde{b}_{12}(\theta)\\\tilde{b}_{21}(\theta)&\tilde{b}_{22}(\theta)\end{pmatrix}.$$
 Then we can rewrite  $B_1(\theta)$ as
$$
B_1(\theta)=\tilde{B}(\theta)R_{\frac{\langle \ell,\theta\rangle}{2}}e^{Y(\theta)}M^{-1}U^{-1}=\tilde{B}(\theta)M^{-1}MR_{\frac{\langle \ell,\theta\rangle}{2}}M^{-1}Me^{Y(\theta)}M^{-1}U^{-1}.
$$
Direct computation shows that
\begin{align*}
z_{11}(\theta)&=e^{-2\pi i\frac{\langle \ell+\ell_0,\theta\rangle}{2}}b_{11}(\theta)\\
&=(u_{11} y_{11}(\theta)+ u_{21}y_{12}(\theta) ) \tilde{b}_{11}(\theta)e^{-2\pi i\frac{\langle \ell_0,\theta\rangle}{2}}\\
&+(u_{11}y_{21}(\theta)+u_{21}y_{22}(\theta))\tilde{b}_{12}(\theta)e^{-2\pi i\frac{\langle \ell_0,\theta\rangle}{2}}e^{-2\pi i\langle \ell,\theta\rangle}.
\end{align*}
Since  $\|Y\|_{h'} \leq  \frac{1}{2},$  then we have
$$\|y_{12}\|_{h'},\|y_{21}\|_{h'}\leq  \|Y\|_{h'}.$$
Consequently, one can estimate $|\hat{z}_{11}(n)|$ as
\begin{eqnarray*}
|\hat{z}_{11}(n)| &\leq&  \|U\|\|\tilde{B}\|_{h'}e^{\pi|\ell_0|h'}e^{-2\pi h'|n|}\\
&+& (|u_{11}| \|y_{21}\|_{h'} +|u_{21}| \|y_{22}\|_{h'} )\|\tilde{B}\|_{h'}e^{\pi|\ell_0|h'}e^{-2\pi h'|n+\ell|}.
\end{eqnarray*}

By $(\ref{f2})$ and the fact that  $|\det(B_1(\theta))|=1$,
one has
$$
2\|b_{11}\|_{L^2}=\|b_{11}\|_{L^2}+\|b_{21}\|_{L^2}\geq\|B_1\|_{C^0}^{-1},
$$
and then we have
$$\|\hat{z}_{11}\|_{\ell^2}=\|\hat{b}_{11}\|_{\ell^2}=\|b_{11}\|_{L^2}\geq(2\|B_1\|_{C^0})^{-1}.
$$
Therefore one can further compute   $|\hat{z}_{11}(n)|$ as
\begin{eqnarray}\nonumber\frac{|\hat{z}_{11}(n)|}{\|\hat{z}_{11}\|_{\ell^2}} &\leq& 2 \|U\|\|\tilde{B}\|_{h'}\|B_1\|_{C^0}e^{\pi|\ell_0|h'}e^{-2\pi  h'|n|}\\
\nonumber&+& 2(|u_{11}|\|y_{21}\|_{h'} +|u_{21}| \|y_{22}\|_{h'} )\|\tilde{B}\|_{h'}\|B_1\|_{C^0}e^{\pi|\ell_0|h'}e^{-2\pi h'|n+\ell|} \\
\label{z11} &\leq& 4 \|U\|^2\|\tilde{B}\|^2_{h'}e^{\pi|\ell_0|h'}e^{-2\pi h'|n|}\\
\nonumber&+& 4(|u_{11}|\|y_{21}\|_{h'} +|u_{21}| \|y_{22}\|_{h'} )\|U\|\|\tilde{B}\|^2_{h'}e^{\pi|\ell_0|h'}e^{-2\pi h'|n+\ell|}.
\end{eqnarray}

Now we distinguish the proof into two cases:

\text{Case 1:} $\frac{2|\nu|}{\|2\rho(E)-\langle \ell,\alpha\rangle-\langle \ell_0,\alpha\rangle\|_{\R/\Z}} \geq \frac{1}{10}$. In this case, since $U$ is unitary, we have
\begin{align*}
\frac{|\hat{z}_{11}(n)|}{\|\hat{z}_{11}\|_{\ell^2}} \leq 8\|\tilde{B}\|_{h'}^{2}e^{\pi|\ell_0|h'}e^{-2\pi h'|n|}+8\|\tilde{B}\|_{h'}^{2}e^{\pi|\ell_0|h'}e^{-2\pi h'|n+\ell|}.
\end{align*}

\text{Case 2:} $\frac{2|\nu|}{\|2\rho(E)-\langle \ell,\alpha\rangle-\langle \ell_0,\alpha\rangle\|_{\R/\Z}} < \frac{1}{10}$. In this case,  we need the following result:

\begin{Lemma}\label{Dominate matrix}
Assume that $$MAM^{-1}=exp \left(
\begin{array}{ccc}
 i t &  \nu\\
\bar{ \nu} &  -i t
 \end{array}\right)\in SU(1,1)$$ with $spec\{A\}=\{e^{2\pi i\rho},e^{-2\pi i\rho}\}$, if $|\frac{4\nu}{\rho}|\leq 1$ and $\rho t>0$, then there exists $U\in SL(2,\C)$, such that
$$
UMAM^{-1}U^{-1}=\begin{pmatrix}e^{2\pi i \rho}&0\\ 0&e^{-2\pi i \rho}\end{pmatrix},
$$
with
$$
\|U-id\|\leq |\frac{\nu}{\rho}|.
$$
\end{Lemma}
\begin{pf} We just need to find some $U\in SL(2,\C)$ such that \begin{eqnarray}\label{domi}U^{-1}\left(
\begin{array}{ccc}
 i t &  \nu\\
\bar{ \nu} &  -i t
 \end{array}\right)U=\begin{pmatrix}i\rho&0\\0&-i\rho\end{pmatrix},\end{eqnarray}
 Note that \eqref{domi} implies $t^2-|\nu|^2=\rho^2$, since $\rho t>0$, we have
$$
|t+\rho|\geq 2|\rho|.
$$
One can easily check that
 $$U=\frac{1}{\sqrt{1-|\frac{\nu}{t+\rho}|^2}}\begin{pmatrix}1&\frac{\nu i}{t+\rho}\\\frac{-\bar{\nu}i}{t+\rho}&1\end{pmatrix}$$ will satisfy our needs.
Moreover,
since $|\frac{1}{\sqrt{1-x^2}}-1|\leq |x|^2$, if $|x|<1$,  then we have
$$
\|U-id\|\leq |\frac{\nu}{t+\rho}|^2+ \frac{3}{2}|\frac{\nu}{t+\rho}|\leq |\frac{\nu}{\rho}|,
$$
where we use the simple fact that $|\frac{\nu}{t+\rho}|\leq \frac{1}{8}$.
\end{pf}

Thus if $(\rho(E)-\frac{\langle \ell+\ell_0,\alpha\rangle}{2})t>0$, by Lemma \ref{Dominate matrix},  instead of taking $U$ to be unitary, one can take $U\in SL(2,\C)$ such that
$$
UMAM^{-1}U^{-1}=\begin{pmatrix}e^{2\pi i(\rho(E)-\frac{\langle \ell+\ell_0,\alpha\rangle}{2})}&0\\0&e^{2\pi i(-\rho(E)+\frac{\langle \ell+\ell_0,\alpha\rangle}{2})}\end{pmatrix},
$$
with estimate $\|U^{-1}-id\|\leq\frac{2|\nu|}{\|2\rho(E)-\langle \ell,\alpha\rangle-\langle \ell_0,\alpha\rangle\|_{\R/\Z}}<\frac{1}{10}$, i.e.
$$|u_{12}|,|u_{21}|\leq \frac{2|\nu|}{\|2\rho(E)-\langle \ell,\alpha\rangle-\langle \ell_0,\alpha\rangle\|_{\R/\Z}}.$$

Similarly, using the basic relation \eqref{z11}, one can estimate
\begin{eqnarray*}
\frac{|\hat{z}_{11}(n)|}{\|\hat{z}_{11}\|_{\ell^2}}&&\leq 8\|\tilde{B}\|_{h'}^{2}e^{\pi|\ell_0|h'}e^{-2\pi h'|n|} \\ &&+ 8\|\tilde{B}\|_{h'}^{2}e^{\pi|\ell_0|h'}(\|Y\|_{h'}+\frac{2|\nu|}{\|2\rho(E)-\langle \ell,\alpha\rangle-\langle \ell_0,\alpha\rangle\|_{\R/\Z}})e^{-2\pi h'|n+\ell|}.
\end{eqnarray*}

Note that
\begin{align*}
\tilde{B}(\theta)&=\sum\limits_{k=-n_0}^{n_0}\hat{\tilde{B}}(k)e^{2\pi i\langle k,\theta\rangle}+\sum\limits_{|k|> n_0}\hat{\tilde{B}}(k)e^{2\pi i\langle k,\theta\rangle}\\
&=T_{n_0}\tilde{B}+R_{n_0}\tilde{B}.
\end{align*}
If $n_0=\frac{\ln\|\tilde{B}\|_{h'}}{\pi h'}$, then
$$
\|\tilde{B}-T_{n_0}\tilde{B}\|_{C^0}=\|R_{n_0}\tilde{B}\|_{C^0} \leq \sum\limits_{|k|>n_0}\|\tilde{B}\|_{h'}e^{-2\pi kh'}\leq \|\tilde{B}\|_{C^0}.
$$
By Rouche's theorem,  this implies that
$$
\deg{\tilde{B}}=\deg{T_{n_0}\tilde{B}},
$$
since the topological degree of $T_{n_0}\tilde{B}$ can be bounded by the numbers of zeros of a non-vanishing coordinate, hence
$$
\deg{T_{n_0}\tilde{B}}\leq 2n_0=\frac{2\ln\|\tilde{B}\|_{h'}}{\pi h'}
$$
which means
$$
e^{\pi|\ell_0|h'}\leq \|\tilde{B}\|_{h'}^2.
$$

Combing this with the above two cases, we obtained that the long-range operator $L_{V,\lambda,\alpha,\rho(E)}$ has a $(2\pi h',\ell, C,C_{|\ell|})$-good eigenfunction.

If $(\rho(E)-\frac{\langle \ell+\ell_0,\alpha\rangle}{2})t<0$, one can take $U'\in SL(2,\C)$ such that
$$
U'MAM^{-1}U'^{-1}=\begin{pmatrix}e^{-2\pi i(\rho(E)-\frac{\langle \ell+\ell_0,\alpha\rangle}{2})}&0\\0&e^{2\pi i(\rho(E)-\frac{\langle \ell+\ell_0,\alpha\rangle}{2})}\end{pmatrix},
$$
with estimate $\|U'^{-1}-id\|\leq\frac{2|\nu|}{\|2\rho(E)-\langle \ell,\alpha\rangle-\langle \ell_0,\alpha\rangle\|_{\R/\Z}}<\frac{1}{10}$.

Let $B'_1(\theta)=\tilde{B}(\theta)R_{\frac{\langle \ell,\theta\rangle}{2}}e^{Y(\theta)}M^{-1}U'^{-1}$, then
$$
\overline{B'_1(\theta+\alpha)}^{-1}S_E^{\lambda^{-1} V}(\theta)\overline{B'_1(\theta)}=\begin{pmatrix}e^{2\pi i(\rho(E)-\frac{\langle \ell+\ell_0,\alpha\rangle}{2})}&\bar{c}\\0&e^{2\pi i(-\rho(E)+\frac{\langle \ell+\ell_0,\alpha\rangle}{2})}\end{pmatrix}.
$$
Since $\overline{B'_1(\theta)}$ has the same structure and estimate as $B_1(\theta)$, by similar argument as above, one can prove that the long-range operator $L_{V,\lambda,\alpha,\rho(E)}$ has a $( 2\pi h',-\ell,C,C_{|\ell|})$-good eigenfunction.
\end{pf}

\section{Quantitative almost reducibility}
In section 4, we have proved that nice reducibility result of Schr\"odinger cocycles implies the dual systems have good eigenfunctions. In this section, we give quantitative reducibility results with desired estimates.
\subsection{Local  quantitative  almost reducibility}

\begin{Proposition}\label{reducibility}
For any $0<\tilde{h}<h$, $\kappa>0,\kappa'>0$, $\tau>d-1$.
Suppose that $\alpha\in DC_d(\kappa',\tau)$, $\rho(\alpha,A_0e^{f_0})\in DC_\alpha(\kappa,\tau)$. Then there exist $B\in C_{\tilde{h}}^\omega(\T^d, PSL(2,\R))$ and $A\in SL(2,\R)$ satisfying
$$
B^{-1}(\theta+\alpha)A_0e^{f_0(\theta)}B(\theta)=A=M^{-1}exp \left(
\begin{array}{ccc}
 i t &  \nu\\
\bar{ \nu} &  -i t
 \end{array}\right)M,
$$  provided that $\|f_0\|_h<\epsilon_*$ for some $\epsilon_*>0$ depending on $A_0,\kappa',\tau,h,\tilde{h},d$.

 In particular, there exist $Y\in (\T^d,s\ell(2,\R))$, $\ell\in \Z^d$, $\tilde{B}\in C_{\tilde{h}}^\omega(\T^d, PSL(2,\R))$, such that $B(\theta)=\tilde{B}(\theta)R_{\frac{\langle \ell,\theta\rangle}{2}}e^{Y(\theta)}$ with the following estimates
\begin{equation}\label{es1}
 \|Y\|_{\tilde{h}}\leq  e^{-2\pi|\ell|\tilde{h}}, \\ \ |\nu|\leq 2e^{-2\pi|\ell|\tilde{h}},
\end{equation}
\begin{equation}\label{es2}
\|\tilde{B}\|_{\tilde{h}}<|\ln\kappa|^{\tau}\kappa^{-\frac{2(h-\tilde{h})}{\tilde{h}}}, \end{equation}
\begin{equation}\label{es3}
 |\ell|\leq \frac{|\ln \kappa|^4}{h-\tilde{h}},
\end{equation}
\begin{equation}\label{eigenvalue}
\|2\rho(\alpha,A_0e^{f_0})-\langle \ell,\alpha\rangle-\langle \deg \tilde{B},\alpha\rangle\|_{\R/\Z}\geq \frac{\kappa}{2^\tau|\ell|^{\tau}}. \end{equation}
\end{Proposition}

\begin{Remark}\label{uniformcons}
 If $A_0$ varies in ${\rm SO}(2,\R)$, then $\epsilon_*$ can be taken uniform with respect to $A_0$.
\end{Remark}
%

\begin{pf} We prove Proposition \ref{reducibility} by iteration. Without loss of generality, we assume that $h<1$. Suppose that
\begin{align*}
\|f_0\|_h\leq \epsilon_* \leq \frac{D_0}{\|A_0\|^{C_0}}(\frac{h-\tilde{h}}{8})^{C_0\tau},\end{align*}
where $D_0=D_0(\kappa',\tau,d)$ is the constant defined in Theorem \ref{iteration}. Then we can define the sequence inductively.  Let
$\epsilon_0=\epsilon_*$, $h_0=h$, assume that we are at the $(j+1)^{th}$ KAM step, i.e. we already construct $B_j\in C^\omega_{h_{j}}(\T^d,PSL(2,\R))$ such that
$$
B_{j}^{-1}(\theta+\alpha)A_0e^{f_0(\theta)}B_{j}(\theta)=A_{j}e^{f_{j}(\theta)},
$$
where $A_j\in SL(2,\R)$ with two eigenvalues $e^{\pm i\xi_j}$ and
$$     \|B_{j}\|_{h_j}\leq \epsilon_j^{-\frac{h-\tilde{h}}{4\tilde{h}}}, \qquad     \|f_j\|_{h_j}\leq \epsilon_j$$
 for some $\epsilon_j\leq \epsilon_0^{2^j}$, then we define
$$
h_j-h_{j+1}=\frac{h-\frac{h+\tilde{h}}{2}}{4^{j+1}}, \ \ N_j=\frac{2|\ln\epsilon_j|}{h_j-h_{j+1}}.
$$

By our selection of  $\epsilon_0$, one can check that
\begin{equation}\label{iter}
\epsilon_j \leq \frac{D_0}{\|A_j\|^{C_0}}(h_j-h_{j+1})^{C_0\tau}.
\end{equation}
 Indeed, $\epsilon_j$ on the left side of the inequality decays at least super-exponentially with $j$, while $(h_j-h_{j+1})^{C_0\tau}$ on the right side decays exponentially with $j$.

Note that $(\ref{iter})$ implies that  Proposition \ref{iteration} can be applied iteratively, consequently one can construct
$$
\bar{B}_j\in C^\omega_{h_{j+1}}(\T^d,PSL(2,\R)),\ \ A_{j+1}\in SL(2,\R),\ \ f_{j+1}\in C_{h_{j+1}}(\T^d,sl(2,\R))
$$
such that
$$
\bar{B}_j^{-1}(\theta+\alpha)A_je^{f_j(\theta)}\bar{B}_j(\theta)=A_{j+1}e^{f_{j+1}(\theta)}.
$$
More precisely, we can distinguish two cases:\\

\noindent \textbf{Non-resonant case:}  If for any $n\in \Z^{d}$ with $0< |n| \leq N_j$, we have
$$
\| 2\xi_j - <n,\alpha> \|_{\R/\Z}\geq \epsilon_j^{\frac{1}{15}},
$$
then
\begin{equation*}
\| \bar{B}_j-id\|_{h_{j+1}}\leq \epsilon_j^{\frac{1}{2}} ,\   \ \| f_{j+1}\|_{h_{j+1}}\leq \epsilon_j^2:=  \epsilon_{j+1}, \ \ \|A_{j+1}-A_j\|\leq 2\epsilon_j.
\end{equation*}
Let $B_{j+1}=B_j(\theta)\bar{B}_j(\theta)$, we have
$$
B_{j+1}^{-1}(\theta+\alpha)A_0e^{f_0(\theta)}B_{j+1}(\theta)=A_{j+1}e^{f_{j+1}(\theta)},
$$
with estimate $$
\|B_{j+1}\|_{h_{j+1}}\leq 2\epsilon_j^{-\frac{h-\tilde{h}}{4\tilde{h}}}\leq \epsilon_{j+1}^{-\frac{h-\tilde{h}}{4\tilde{h}}}.
$$
Moreover,
\begin{equation}\label{deg1}
\deg{B_{j+1}}=\deg{B_{j}},\end{equation}
 since $\bar{B}_j(\theta)$ is close to the identity. \\
%

\noindent \textbf{Resonant case:} If there exists $n_j$ \footnote{We call such $n_j$ the resonance.} with $0<| n_j| \leq N_j$ such that
$$
\| 2\xi_j- <n_j,\alpha> \|_{\R/\Z}< \epsilon_j^{\frac{1}{15}},
$$
then $\bar{B}_j(\theta)=\bar{B}_j'(\theta)R_{\frac{\langle n_j,\theta\rangle}{2}}$ with estimate
$$\|\bar{B}_j\|_{h_{j+1}}\leq |\ln\epsilon_j|^{\tau}\epsilon_j^{-\frac{h_{j+1}}{h_j-h_{j+1}}}, \ \  \|\bar{B}'_j\|_{h_{j+1}}< |\ln\epsilon_j|^{\tau}.$$
$$ \| f_{j+1}\|_{h_{j+1}}  \leq \epsilon_j e^{-h_{j+1}\varepsilon_j^{-\frac{1}{18\tau}}} := \epsilon_{j+1}.$$
Moreover, we can write

 $$A_{j+1}= M^{-1}exp \left(
\begin{array}{ccc}
 i t^{j+1} &  \nu^{j+1}\\
\bar{ \nu}^{j+1} &  -i t^{j+1}
 \end{array}\right)M
 $$ with estimate
$$ |\nu^{j+1}|< e^{-2\pi|n_j|h_j}.$$
Let $B_{j+1}(\theta)=B_j(\theta)\bar{B}_j(\theta)$, then we have
$$
B_{j+1}^{-1}(\theta+\alpha)A_0e^{f_0(\theta)}B_{j+1}(\theta)=A_{j+1}e^{f_{j+1}(\theta)},
$$
with
\begin{equation}\label{deg2}
\deg{B_{j+1}} =\deg{B_{j}}+\deg \bar{B}_j= \deg{B_{j}}  + n_j
\end{equation}
\begin{align*}
\|B_{j+1}\|_{h_{j+1}}&\leq \epsilon_j^{-\frac{h-\tilde{h}}{4\tilde{h}}}|\ln\epsilon_j|^\tau\epsilon_j^{-\frac{h_{j+1}}{h_j-h_{j+1}}}
\leq \epsilon_{j+1}^{-\frac{h-\tilde{h}}{4\tilde{h}}}.
\end{align*}
The last inequality is possible since by our selection $ \epsilon_{j+1}=\epsilon_j e^{-h_{j+1}\varepsilon_j^{-\frac{1}{18\tau}}}$.

Let $\ell_j=\deg{B_j}$, by \eqref{deg1} and \eqref{deg2}, one has
\begin{equation}\label{degree}
|\ell_{j}|=|\deg{B_{j}}|\leq \sum\limits_{i=0}^{j-1}N_i\leq \sum\limits_{i=0}^{j-1}\frac{2|\log\epsilon_{i}|}{h_i-h_{i+1}}\leq \frac{ |\log\epsilon_{j-1}|^3}{h-\tilde{h}}.
\end{equation}

On the other hand,  if $\phi\in DC_{\alpha}(\kappa,\tau)$,   then  by \eqref{degree},
we have
\begin{eqnarray*}
&&\|2\phi-\langle m,\alpha\rangle-\langle \ell_{j},\alpha\rangle\|_{\R/\Z} \\
&\geq& \frac{\kappa}{(|m+\ell_{j}|+1)^{\tau}}\geq \frac{(1+|\ell_{j}|)^{-\tau}\kappa}{(|m|+1)^{\tau}} \\
& \geq &  \frac{1}{(|m|+1)^{\tau}}. \frac{\kappa (h-\tilde{h})^\tau}{  |\log\epsilon_{j-1}|^{3\tau}},
\end{eqnarray*}
 this implies that
\begin{align}\label{rotation}
\rho(\alpha,A_{j}e^{f_{j}})\in DC_\alpha( \frac{\kappa (h-\tilde{h})^\tau}{  |\log\epsilon_{j-1}|^{3\tau}},\tau).
\end{align}

Now we select $j_0\in\Z$  to be the smallest integer $j$ satisfying
\begin{align}\label{reducibility condition}
\epsilon_{j}\leq \epsilon(\tau,\kappa', \frac{\kappa (h-\tilde{h})^\tau}{  |\log\epsilon_{j-1}|^{3\tau}},\frac{h+\tilde{h}}{2},\tilde{h},d,A_j),
\end{align}
where $\epsilon=\epsilon(\tau,\kappa',\kappa,h,\tilde{h},d,R)$ is defined in Theorem \ref{positive reducibility}.
 Note that in our iteration process, the constant matrix $A_j$ is uniformly bounded, by the definition of $\epsilon$, one can thus select $\epsilon(\tau,\kappa, \frac{\kappa (h-\tilde{h})^\tau}{  |\log\epsilon_{j-1}|^{3\tau}},\frac{h+\tilde{h}}{2},\tilde{h},d,A_j)$ to be independent of $A_j$. Therefore, by the definition of $\epsilon$ and our selection, it  follows that
$$\epsilon_{j_0-1}\geq D_0(\frac{h- \tilde{h}}{2})^{C_0\tau}(\frac{\kappa (h-\tilde{h})^\tau}{  |\log\epsilon_{j_0-2}|^{3\tau}})^{4}.$$
By the fact  $\epsilon_j\leq \epsilon_0^{2^j}$ and the definition of $\epsilon_0$, one therefore has
\begin{equation}\label{estimate}\epsilon_{j_0-1}^{\frac{1}{2}}\geq \kappa^4,
\end{equation}

 Note that \eqref{reducibility condition} means
the condition in Theorem \ref{positive reducibility} is satisfied,  consequently, there exists $\tilde{Y}\in C_{\tilde{h}}^\omega(\T^d,s\ell(2,\R))$ such that
$$
e^{-\tilde{Y}(\theta+\alpha)}A_{j_0}e^{f_{j_0}(\theta)}e^{\tilde{Y}(\theta)}=A=M^{-1}exp \left(
\begin{array}{ccc}
 i t &  \nu\\
\bar{ \nu} &  -i t
 \end{array}\right)M,
$$
with estimates $\|\tilde{Y}\|_{\tilde{h}}\leq \epsilon_{j_0}^{\frac{1}{2}}$, $\|A-A_{j_0}\|\leq \epsilon_{j_0}$.
Let $B=B_{j_0}e^{\tilde{Y}}$, we have
\begin{equation}\label{11}
B^{-1}(\theta+\alpha)A_0e^{f_0(\theta)}B(\theta)=A=M^{-1}exp \left(
\begin{array}{ccc}
 i t &  \nu\\
\bar{ \nu} &  -i t
 \end{array}\right)M.
\end{equation}

Now we prove that the conjugacy $B(\theta)$ can be written in the desired form with good estimates. Indeed, let $\ell$ be the last resonance, and we may assume that the last resonance happens at some step $0<j_0'<j_0$ ( if $j_0'=j_0$ then the proof would be much simpler). By the above iteration process, there exist $\{\bar{B}_j\}_{j=j_0'-1}^{j_0-1}$ and $\{A_j\}_{j=j_0'-1}^{j_0-1}$ such that
\begin{equation}
\bar{B}^{-1}_{j}(\theta+\alpha)A_je^{f_j(\theta)}\bar{B}_{j}(\theta)=A_{j+1}e^{f_{j+1}(\theta)}.
\end{equation}

If $j=j_0'-1$, then by the selection of $j_0'$ and the iteration process, one has  $\bar{B}_{j_0'-1}=\bar{B}_{j_0'-1}'R_{\frac{\langle \ell,\theta\rangle}{2}}$ with
$$\| \bar{B}'_{j_0'-1} \|_{h_{j_0'}}\leq |\ln\epsilon_{j_0'-1}|^{\tau}, \qquad \| f_{j_0'}\|_{h_{j_0'}}\leq  \epsilon_{j_0'}   \ll e^{-2\pi|\ell|\tilde{h}}.$$
Furthermore, we have
 $$MA_{j_0'}M^{-1}= exp \left(
\begin{array}{ccc}
 i t^{j_0'} &  \nu^{j_0'}\\
\bar{ \nu}^{j_0'} &  -i t^{j_0'}
 \end{array}\right)
 $$ with estimate
$|\nu^{j'_0}|< e^{-2\pi|\ell|\tilde{h}}.$

For $j_0'\leq j\leq j_0-1$, by the iteration process, we have
$$
\|\bar{B}_j-id\|_{h_{j+1}}\leq \epsilon_j^{\frac{1}{2}}, \ \ \|A_j-A_{j+1}\|\leq 2\epsilon_j.
$$
Let $\tilde{B}(\theta)=B_{j_0'-1}(\theta)\bar{B}'_{j_0'-1}(\theta)$, $e^{Y(\theta)}=\bar{B}_{j_0'}(\theta)\cdots \bar{B}_{j_0-1}(\theta)e^{\tilde{Y}(\theta)}$, then $B(\theta)$ can be written as $B(\theta)=\tilde{B}(\theta)R_{\frac{\langle \ell,\theta\rangle}{2}}e^{Y(\theta)}$ with estimates
\begin{equation*}
\|Y\|_{\tilde{h}}\leq C\epsilon_{j_0'}^{\frac{1}{2}}\ll e^{-2\pi|\ell|\tilde{h}},
\end{equation*}
$$
 \|A-A_{j_0'}\|\leq \|A-A_{j_0}\|+\|A_{j_0}-A_{j_0'}\|\leq C\epsilon_{j_0'}^{\frac{1}{2}}.
$$
Therefore, it concludes that
\begin{eqnarray*}
|\nu| &\leq& |\nu^{j_0'}| +\|A-A_{j_0'}\| \\
&\leq&e^{-2\pi|\ell|\tilde{h}}+C\epsilon_{j_0'}^{\frac{1}{2}}<2e^{-2\pi|\ell|\tilde{h}},
\end{eqnarray*}
and then \eqref{es1} is proved.

Moreover,  by \eqref{estimate}, one has
\begin{equation*}
\|\tilde{B}\|_{\tilde{h}}<\epsilon_{j_0'-1}^{-\frac{h-\tilde{h}}{4\tilde{h}}}|\ln\epsilon_{j_0'-1}|^{\tau} \leq |\ln\kappa|^{\tau}\kappa^{-\frac{2(h-\tilde{h})}{\tilde{h}}}, \end{equation*}
which  proves \eqref{es2}.

To estimate the remained inequalities, we need more detailed analysis on the resonances. Assume that there are at least two resonant steps, say the $(m_i+1)^{th}$ and $(m_{i+1}+1)^{th}$. At the $(m_{i+1}+1)^{th}$-step, the resonance condition implies $|\xi_{m_{i+1}}-\frac{\langle n_{m_{i+1}},\alpha\rangle}{2}|\leq \epsilon_{m_{i+1}}^{\frac{1}{15}}$, hence $|\xi_{m_{i+1}}|\geq \frac{\kappa}{3|n_{m_{i+1}}|^{\tau}}$. On the other hand, according to Proposition \ref{iteration}, after the $(m_i+1)^{th}$-step, $|\xi_{m_i+1}|\leq \epsilon_{m_i}^{\frac{1}{16}}$. Thus
\begin{equation}\label{different resonances}
|n_{m_{i+1}}|\geq \epsilon_{m_i}^{-\frac{1}{18\tau}}|n_{m_i}|.
\end{equation}

Assuming that there are $s+1$ resonant steps, associated with integers vectors
$$
n_{m_0},...,n_{m_s}=n_{j_0'}=\ell \in\Z^d, \ \ 0<|n_{m_i}|\leq N_{m_i},\ \ i=0,1,...,s,
$$
in view of \eqref{different resonances}, we have
$$
|\deg{\tilde{B}}| =| n_{m_0}+\cdots + n_{m_{s-1}}  |   \leq 2\epsilon_{0}^{\frac{1}{18\tau}}|\ell|.
$$
thus we have $|\ell+\deg{\tilde{B}}|\leq  2|\ell|$, then
\begin{equation*}
\|\rho(\alpha,A_0e^{f_0})-\langle \ell,\alpha\rangle-\langle \deg{\tilde{B}} ,\alpha\rangle\|_{\R/\Z}\geq\frac{\kappa}{|\ell+\deg{\tilde{B}}|^{\tau}}\geq\frac{\kappa}{2^\tau |\ell|^{\tau}}.
\end{equation*}
By the definition of $\ell$ and \eqref{estimate}, one has
$$  | \ell| \leq \frac{2|\ln\epsilon_{j_0'-1}|}{h_{j_0'-1}-h_{j_0'}} \leq \frac{|\ln\kappa|^4}{h-\tilde{h}}.$$
Thus we have finished the whole proof. \end{pf}

\subsection{Global to local reduction}

Proposition \ref{reducibility} deals with reducibility of cocycles close to constant. If $d=1$, one can indeed deal with all subcritical cocycles, with the help of Avila's solution of Almost Reducibility Conjecture \cite{avila1,avila2}.

\begin{Lemma}\label{Global to local general}
Let  $\alpha\in {\rm DC}_1$, $V\in C^{\omega}(\T,\R)$ with $\Sigma_{V,\alpha}=\Sigma_{V,\alpha}^{sub}$.  There exist $h_1=h_1(\alpha,V)>0$ such that for any $\eta>0$, $E\in \Sigma_{V,\alpha}$, there exist $\Phi_{E}\in C^{\omega}(\T, PSL(2,\R))$ with $|\Phi_{E}|_{h_1}<\Gamma(V,\alpha,\eta)$ such that
\begin{equation}
\Phi_{E}(\theta+\alpha)^{-1}S_E^{V}(\theta)\Phi_{E}(\theta)=R_{\phi(E)}e^{f_E(\theta)},
\end{equation}
with $\|f_E\|_{h_1}<\eta$.
\end{Lemma}

\begin{pf}
The crucial fact in this proposition is that we can choose $h_1$ to be independent of $E$ and $\eta$, and choose $\eta$ to be independent of $E$.  The ideas of the proof are essentially contained in Proposition 5.1 and Proposition 5.2 of \cite{LYZZ}, we include the proof here for completeness.

For any $E\in\Sigma^{\rm sub}_{V,\alpha}$, the cocycle $(\alpha, S_{E}^V)$ is subcritical, hence almost reducible by Theorem  \ref{arc-conjecture}, thus there exists  $h_0= h_0(E,V,\alpha)>0$, such that  for any $\eta>0$, there are $\Phi_E\in C_{h_0}^\omega(\T,{\rm PSL}(2,\R))$, $F_E\in C_{h_0}^\omega(\T,{\rm gl}(2,\R))$  and $\phi(E)\in\T$ such that
$$\Phi_E(\cdot+\alpha)^{-1} S_E^{V}(\cdot) \Phi_E(\cdot)=R_{\phi(E)}+F_E(\cdot),$$
with $\|F_E\|_{h_0}< \eta/2$ and $\|\Phi_E\|_{h_0}<\tilde{\Gamma
}$ for some $\tilde{\Gamma}=\Gamma(V,\alpha, \eta, E)>0$.
As a consequence, for any $E'\in \R$, one has
$$
\left\|\Phi_E(\cdot+\alpha)^{-1} S_{E'}^{V}(\cdot) \Phi_E(\cdot) -  R_{\phi(E)}\right\|_{h_0} < \frac{\eta}{2}+ |E-E'| \, \| \Phi_E\|_{h_0}^2.
$$
It follows that with the same $\Phi_E$, we have $$\|\Phi_E(\theta+\alpha)^{-1} S_{E'}^{V}(\theta) \Phi_E(\theta) - R_{\phi(E)}\|_{h_0} < \eta$$ for any energy $E'$ in a neighborhood $\mathcal{U}(E)$ of $E$.  Since $\Sigma_{V,\alpha}$ is compact, by compactness argument, we can select  $h_0(E,V,\alpha)$, $\Gamma(V,\alpha, \eta, E)>0$ to be   independent of the energy $E$.
\end{pf}

Note that for general subcritical Schr\"odinger cocycles, one can only obtain the existence of $h_1$, however for almost Mathieu cocycles, one can give a good control of $h_1$. This is the main reason why we can give sharp exponential decay rate in expectation for  almost Mathieu operators.

\begin{Lemma}[\cite{LYZZ}]\label{amoredu}
Let  $\alpha\in {\rm DC}_1$, $\lambda>1$.  For any $\e>0$,  $\eta>0$, if $E\in \Sigma_{\lambda^{-1},\alpha}$, there exist $\Phi_{E}\in C^{\omega}(\T, PSL(2,\R))$ with $|\Phi_{E}|_{\frac{1}{2\pi}\ln\lambda-\e}<\Gamma$ for some $\Gamma=\Gamma(\lambda,\alpha,\eta,\e)>0$ such that
\begin{equation}
\Phi_{E}(\theta+\alpha)^{-1}S_E^{2\lambda^{-1}\cos}(\theta)\Phi_{E}(\theta)=R_{\phi(E)}e^{f_E(\theta)},
\end{equation}
with $\|f_E\|_{\frac{1}{2\pi}\ln\lambda-\e}<\eta$.
\end{Lemma}

\section{More precise estimates on good eigenfunctions}
In this section, we give more precise estimates on $(\gamma, \ell, C,C_{|l|})$  by  the quantitative almost reducibility estimates   given in section 5.

\begin{Proposition}\label{13}
Let $\alpha\in DC_1$, $V\in C^{\omega}(\T,\R)$ with $\Sigma_{\lambda^{-1} V,\alpha}=\Sigma_{\lambda^{-1} V,\alpha}^{sub}$. Then for any $\e>0$, there exist $h_1=h_1(V,\alpha)$, $C_4=C_4(V,\alpha,\e)$ with the following properties: if
$\rho(E(\omega))=\omega\in DC_\alpha(\kappa,\tau)$,  then associated with the eigenvalue $\lambda E(\omega)$,  the long range operator $L_{V,\lambda,\alpha,\omega}$ has a \small $$(2\pi (h_1- \frac{\e}{96}),\ell,C_4|\ln\kappa|^{4\tau}\kappa^{-\frac{\e}{10 h_1}},\min\{1,\frac{C_4 |\ell|^{\tau}e^{-2\pi|\ell|(h_1-\frac{\e}{96})}}{\kappa}\})$$-good eigenfunction for some $|\ell|\leq C_4|\ln\kappa|^4$.  \end{Proposition}

\begin{pf}
Since  $h_1(\alpha,V)$ in Lemma \ref{Global to local general} is fixed, and it is independent of $\eta$, then  one can always take $\eta$ small enough such that
$$\eta \leq \epsilon_*(\kappa',\tau,h_1,h_1-\e/96,d),$$
where $\epsilon_*(A_0,\kappa',\tau,h,\tilde{h},d)$ is the constant defined in  Proposition \ref{reducibility}.
Note that by Remark \ref{uniformcons}, the constant
$\epsilon_*$ given by Proposition \ref{reducibility} can be taken uniformly with respect to $R_{\phi} \in {\rm SO}(2,\R)$.
By Lemma \ref{Global to local general}, there exists $\Phi_{E(\omega)}\in C_{h_1}^\omega(\T, PSL(2,\R))$ such that
$$\Phi_{E(\omega)}^{-1}(\theta+\alpha)S_{E(\omega)}^{\lambda^{-1} V}(\theta)\Phi_{E(\omega)}(\theta)=R_{\phi(E(\omega))} e^{f_{E(\omega)}(\theta)}.$$
By footnote 5 of \cite{avila1}, $|\deg{\Phi_{E(\omega)}}|\leq  C |\ln \Gamma|:=\Gamma_1$ for some constant $C=C(V,\alpha)>0$.
Since $\rho(E(\omega))=\omega\in DC_{\alpha}(\kappa,\tau)$, by similar argument as in \eqref{rotation}, one has
$$
\rho(\alpha,\Phi_{E(\omega)}^{-1}(\theta+\alpha)S_{E(\omega)}^{\lambda^{-1} V}(\theta)\Phi_{E(\omega)}(\theta))\in DC_\alpha( \frac{\kappa}{\Gamma_1^{\tau}},\tau).
$$

By our selection, $\|f_{E(\omega)}\|_{h_1}\leq \eta \leq \epsilon_*(\kappa',\tau,h_1,h_1-\e/96,d)$, then we can apply Proposition \ref{reducibility}, and obtain $B\in C_{h_1-\frac{\e}{96}}^\omega(\T, PSL(2,\R))$ and $A\in SL(2,\R)$, such that
$$
B(\theta+\alpha)^{-1}R_{\phi(E(\omega))}e^{f_{E(\omega)}(\theta)}B(\theta)=A=M^{-1}exp \left(
\begin{array}{ccc}
 i t &  \nu\\
\bar{ \nu} &  -i t
 \end{array}\right)M,
$$
Moreover, $B(\theta)$ can be written as $B(\theta)=\tilde{B}(\theta)R_{\frac{\langle \ell,\theta\rangle}{2}}e^{Y(\theta)}$ with the following estimates
\begin{equation*}
 \|Y\|_{h_1-\frac{\e}{96}}\leq  e^{-2\pi|\ell|(h_1-\frac{\e}{96})}, \\ \ |\nu|\leq 2e^{-2\pi|\ell| (h_1-\frac{\e}{96})},
\end{equation*}
\begin{equation*}
\|\tilde{B}\|_{h_1-\frac{\e}{96}}<|\ln\kappa\Gamma_1^{-\tau}|^{\tau}(\kappa \Gamma_1^{-\tau})^{-\frac{\e}{40 h_1}},  \end{equation*}
\begin{equation*}
 |\ell|\leq \frac{96|\ln \kappa\Gamma_1|^4}{\e} \leq   C_4(V,\alpha,\e)|\ln\kappa|^4 ,\end{equation*}
\begin{equation}\label{es5}
\|2\rho(\alpha,R_{\phi(E(\omega))}e^{f_{E(\omega)}})-\langle \ell,\alpha\rangle-\langle \deg \tilde{B},\alpha\rangle\|_{\R/\Z}\geq \frac{\kappa}{2^\tau\Gamma_1^\tau|\ell|^{\tau}}. \end{equation}

Let $B_2(\theta)=B_1(\theta)R_{\frac{\langle \ell,\theta\rangle}{2}}e^{Y(\theta)}=\Phi_{E(\omega)}(\theta)\tilde{B}(\theta)R_{\frac{\langle \ell,\theta\rangle}{2}}e^{Y(\theta)}$, then
$$
B_2(\theta+\alpha)^{-1}S_{E(\omega)}^{\lambda^{-1} V}(\theta)B_2(\theta)=A=M^{-1}exp\begin{pmatrix}it&\nu\\ \bar{\nu}&-it\end{pmatrix}M,
$$
with estimate
\begin{eqnarray*} && \|B_1\|^4_{h_1-\frac{\e}{96}} = \|\Phi_{E(\omega)}\tilde{B}\|^4_{h_1-\frac{\e}{96}}\\  & \leq& \Gamma^4  |\ln\kappa\Gamma_1^{-\tau}|^{4\tau}(\kappa\Gamma_1^{-\tau})^{-\frac{\e}{ 10 h_1}}  \leq C_4(V,\alpha,\e) |\ln\kappa|^{4\tau}  \kappa^{-\frac{\e}{10 h_1}}.
\end{eqnarray*}
Moreover, \eqref{es5} can be written as  \begin{equation*}
\|2\rho(E(\omega))-\langle \ell,\alpha\rangle-\langle \deg{B_1},\alpha\rangle\|_{\R/\Z}\geq \frac{\kappa}{2^\tau\Gamma_1^{\tau}  |\ell|^{\tau}},
\end{equation*}
and consequently
\begin{align*}
&\ \ \ \ \|Y\|_{h_1-\frac{\epsilon}{96}}+\frac{2|\nu|}{\|2\rho(E(\omega))-\langle k,\alpha\rangle-\langle \deg{B_1},\alpha\rangle\|_{\R/\Z}}\\
&\leq e^{-2\pi|\ell|(h_1-\frac{\e}{96})}+  \frac{ 2^{2+\tau}\Gamma_1^{\tau}  |\ell|^{\tau}   e^{-2\pi|\ell|(h_1-\frac{\e}{96})}}{\kappa}\\
&\leq \frac{C_4(V,\alpha,\e) |\ell|^{\tau} e^{-2\pi|\ell|(h_1-\frac{\e}{96})}}{\kappa}.
\end{align*}

By Proposition \ref{Upper bound}, the long range  operator $L_{V,\lambda,\alpha,\omega}$ has a \small$$(2\pi (h_1- \frac{\e}{96}),\ell, C_4|\ln\kappa|^{4\tau}\kappa^{-\frac{\e}{10 h_1}},\min\{1,\frac{C_4 |\ell|^{\tau}e^{-2\pi|\ell|(h_1-\frac{\e}{96})}}{\kappa}\})$$-good eigenfunction for some $|\ell|\leq C_4|\ln\kappa|^4$.
\end{pf}

\begin{Proposition}\label{Amo}
Let $\alpha\in DC_1$, $\lambda>1$.  Then for any $\e>0$, there exists  $C_5=C_5(\lambda,\alpha,\e)$ with the following properties: if
$\rho(E(\omega))=\omega\in DC_\alpha(\kappa,\tau)$, then  associated with the eigenvalue $\lambda E(\omega)$, the almost Mathieu operator $H_{\lambda,\alpha,\omega}$ has a \small $$(2\pi(h_1-\frac{\e}{96}),\ell,C_5|\ln\kappa|^{4\tau}\kappa^{-\frac{\e}{10 h_1}},\min\{1,\frac{C_5 |\ell|^{\tau}e^{-2\pi|\ell|(h_1-\frac{\e}{96})}}{\kappa}\})$$-good eigenfunction for some $|\ell|\leq C_5|\ln\kappa|^4$ with $h_1=\frac{1}{2\pi}\ln\lambda-\e$.\end{Proposition}

\begin{pf}The proof is exactly the same as Proposition \ref{13}. One only need to replace Lemma \ref{Global to local general} by Lemma \ref{amoredu}.
Notice that by Lemma \ref{amoredu}, one can  actually take $h_1=\frac{1}{2\pi}\ln\lambda-\e$, the rest proofs are exactly the same, we omit the details.
\end{pf}

\begin{Proposition}\label{long2}
Let $\alpha\in DC_d$, $V\in C_h^{\omega}(\T,\R)$. Then for any $\e>0$, there exist $\lambda_0(\alpha,V,d,\e)$, $C_6=C_6(V,\alpha,\e)$ with the following properties: if $\lambda>\lambda_0$,
$\rho(E(\omega))=\omega\in DC_\alpha(\kappa,\tau)$, then associated with the eigenvalue $\lambda E(\omega)$, the quasi-periodic long-range operator $L_{V,\lambda,\alpha,\omega}$ has a \small $$(2\pi (h- \frac{\e}{96}),\ell,C_6|\ln\kappa|^{4\tau}\kappa^{-\frac{\e}{10 h}},\min\{1,\frac{C_6 |\ell|^{\tau}e^{-2\pi|\ell|(h-\frac{\e}{96})}}{\kappa}\})$$-good eigenfunction for some $|\ell|\leq C_6|\ln\kappa|^4$.
\end{Proposition}

\begin{pf} In this case, one in fact don't need to take the global to local reduction procedure (Lemma \ref{Global to local general}), or one can say $\Phi_E=\id$ in this case,   thus one can take $h_1=h$, and  the rest of the proof follows exactly the same as the proof of  Proposition \ref{13}.
\end{pf}

\section{Proof of main results}

We  give the proof of Theorem \ref{amo} in detail. The proof of Theorem \ref{long} and \ref{highamo} are almost the same. Before proving Theorem \ref{amo}, we first prove the following more general result.

\begin{Theorem}\label{subedl}
Let $\alpha\in DC_1$, $V\in C^{\omega}(\T,\R)$ with $\Sigma_{\lambda^{-1} V,\alpha}=\Sigma_{\lambda^{-1} V,\alpha}^{sub}$.  Then
$L_{V,\lambda,\alpha,\theta}$ has EDL with exponential decay rate in expectation $\gamma(L) \geq 2\pi h_1$, where $h_1$ is defined in Proposition \ref{13}.
\end{Theorem}
\begin{pf}

First we verify that $(L_{V,\lambda,\alpha,\theta})_{\theta\in\T}$ forms  an ergodic family, since $S_n\theta=\theta+n\alpha$ is ergodic. Now we fix
$\gamma=2\pi(h_1-\frac{\e}{2})$, for $i\geq 0$, we define the sequence
$$ c_i =10^{-i},  \quad \Theta_{c_i}=\{\phi|\phi\in DC_\alpha(c_i,\tau)\}. $$
For any  $\theta \in \cup_i\Theta_{c_i}\cap [0,\frac{1}{2}]$,
there exists exactly one $E$ such that  $\rho(E)=\theta$,
since it is well-known that $\rho$ is nonincreasing and can be constant only on intervals contained in $\R \backslash \Sigma$ on which the rotation number must be rationally dependent with respect to $\alpha$ \cite{H,JM}. Let us denote this inverse function by $E(\theta)$, initially defined on $ \cup_i\Theta_{c_i}\cap [0,\frac{1}{2}]$,  extend it evenly onto $ \cup_i\Theta_{c_i}\cap [-\frac{1}{2},0)$,
and then 1-periodically onto $\cup_i\Theta_{c_i}$.

Now we define the sequence
\begin{eqnarray*}
n_i &=&C_4|\log c_i|^4, \\
C_i&=&C_4|\ln c_i|^{4\tau}c_i^{-\frac{\e}{10h_1}}, \\
C_{i,m}&=&\min\{1,\frac{C_4 m^{\tau}e^{-2\pi m(h_1-\frac{\e}{96})}}{c_i}\}.
\end{eqnarray*}
For any $\theta \in \Theta_{c_i}$, by Proposition \ref{13}, we obtain that  the long range operator $L_{V,\lambda,\alpha,\theta}$ has an eigenfunction
$u(\theta)$ which is  $(\gamma,\ell,C_i,C_{i,|\ell|})$-good for some $|\ell|\leq n_i$.

Note by our construction of $E(\theta)$, for any $\theta\in\cup_i\Theta_{c_i}$, $E(S_m\theta)\neq E(S_n\theta)$ for $m\neq n$,
thus $T_{-m}u(S_m\theta)$, $T_{-n}u(S_n\theta)$ are orthonormal for $m\neq n$.  Hence $(L_{V,\lambda,\alpha,\theta})_{\theta\in\T}$ is an ergodic covariant family and assumption (1) in Theorem \ref{2} is verified.

By Remark \ref{extra}, $(L_{V,\lambda,\alpha,\theta})_{\theta\in\T}$ display AL for a.e. $\theta$. Since $E(S_m\theta)\neq E(S_n\theta)$ for $m\neq n$,  we furthermore obtain that  the point spectrum is simple. Since $E(\theta)$ is measurable, thus $u(\theta)$ can be chosen  to be a measurable function. One can consult \cite{jk2} for more details about this fact. This means assumption (2) in Theorem \ref{2} is satisfied.

Now we  verify the assumption $(3)$ in Theorem \ref{2}. Note that there exists $C(\e)$ such that $m^\tau<e^{2\pi m\e/96}$ for $m>C(\e)$. Thus if we let $\tilde{C}_4=C_4C(\e)^\tau$, then $$C_4 m^{\tau}e^{-2\pi m(h_1-\frac{\e}{96})}\leq \tilde{C}_4 e^{-2\pi m(h_1-\frac{\e}{48})},$$ hence we have
\begin{align*}
\sup\limits_{m}C_{i,m}e^{\gamma m}&\leq\sup\limits_{m\leq \frac{ |\ln c_i /\tilde{C}_4|}{2\pi(h_1-\e/48)}}C_{i,m}e^{\gamma m}+\sup\limits_{m\geq\frac{ |\ln c_i /\tilde{C}_4|}{2\pi(h_1-\e/48)}}C_{i,m}e^{\gamma m}\\
&\leq \sup\limits_{m\leq \frac{ |\ln c_i /\tilde{C}_4|}{2\pi(h_1-\e/48)}}e^{\gamma m}+\sup\limits_{m\geq\frac{ |\ln c_i /\tilde{C}_4|}{2\pi(h_1-\e/48)}}\frac{\tilde{C}_4e^{-2\pi m(h_1-\frac{\e}{48})}}{c_i}e^{\gamma m}\\
&\leq (\frac{\tilde{C}_4}{c_i})^{1-\frac{2\e}{5h_1}}+(\frac{\tilde{C}_4}{c_i})^{1-\frac{2\e}{5h_1}}=2(\frac{\tilde{C}_4}{c_i})^{1-\frac{2\e}{5h_1}}.\\
\end{align*}
It follows  that \begin{align*}
&\ \ \ \ \sum\limits_{i=0}^\infty C_i^2(1+\sup\limits_{m\leq n_i}C_{i,m}e^{m\gamma})c_{i-1}\leq \sum\limits_{i=0}^\infty C_i^2(1+2(\frac{\tilde{C}_4}{c_i})^{1-\frac{2\e}{5h_1}})c_{i-1}\\
&\leq \sum\limits_{i=0}^\infty 10C_4^2|\ln c_i|^{8\tau}c_i^{-\frac{\e}{5h_1}}c_{i}+20\sum\limits_{i=0}^\infty C_4^2|\ln c_i|^{8\tau}c_i^{-\frac{\e}{5h_1}}(\frac{\tilde{C}_4}{c_i})^{1-\frac{2\e}{5h_1}}c_{i}\\
&\leq \sum\limits_{i=0}^\infty 10C^2_4|\ln c_i|^{8\tau}c_i^{1-\frac{\e}{5h_1}}+20\sum\limits_{i=0}^\infty \tilde{C}_4C^2_4|\ln c_i|^{8\tau}c_i^{\frac{\e}{5h_1}}<\infty.
\end{align*}
Thus by Theorem \ref{2}, one obtain that  $\gamma(H)\geq 2\pi(h_1-\frac{\e}{2})$.

We remark that   Propostion \ref{13} holds for any fixed $\e>0$, since one can always make a global to local reduction by Lemma \ref{Global to local general} no matter how small $\e$ is, and then use the local estimates in Proposition \ref{reducibility}. By the definition of exponential decay rate in expectation, one therefore obtain $\gamma(H)\geq 2\pi h_1$.

\end{pf}

\begin{Corollary}\label{nonperturbativelong}
If $\alpha\in DC_1$,   $V\in C^{\omega}(\T,\R)$. Then for any $\e>0$, there exists $\lambda_0(V,\e)$, such that if $\lambda>\lambda_0$,
$L_{V,\lambda,\alpha,\theta}$ has EDL.
\end{Corollary}
\begin{pf}
By Theorem 3.5 in \cite{aj1}, there exists $\lambda_0(V,\e)$, such that $(\alpha.S_E^{\lambda^{-1} V})$ is subcritical if $\lambda>\lambda_0$. One can also prove this by upper semi-continuity of the acceleration \cite{avila0}. Then the result follows from Theorem \ref{subedl} directly.
\end{pf}
\textbf{Proof of Theorem \ref{amo}}. Replace Proposition \ref{13}  by Proposition \ref{Amo},  using same argument as in Theorem \ref{subedl}, one can prove that   for $\alpha\in DC_1$, $\lambda>1$, then $H_{\lambda,\alpha,\theta}$ has EDL  with exponential decay rate in expectation $\gamma(H)\geq \ln\lambda$. Thus we only need to prove $\gamma(H)\leq \ln\lambda$.
We prove it by contradiction. Suppose that there exists $\e_0>0$ and $N_0(\e_0)$ such that for $n>N_0$, we have
\begin{equation}\label{7.2}
\int_{\Omega}\sup\limits_{t\in\R}|\langle\delta_{n}, e^{-itH_\omega}\delta_{0}\rangle|d\omega \leq e^{-(\ln\lambda+\e_0)|n|}.
\end{equation}
Let
$$
Q(\omega)=\sum\limits_{n>N_0}e^{(\ln\lambda+\e_0/2)|n|}\sup\limits_{t\in\R}|\langle\delta_{n}, e^{-itH_\omega}\delta_{0}\rangle|.
$$
Then  \eqref{7.2} implies that
$$
\int Q(\omega) d\omega<\infty.
$$
It follows that  $Q(\omega)<\infty$ for a.e. $\omega$,
$$
\sup\limits_{t\in\R}|\langle\delta_{n}, e^{-itH_\omega}\delta_{0}\rangle|\leq C_\omega e^{-(\ln\lambda+\frac{\e_0}{2})|n|},
$$
provide $n>N_0$.

This implies that that $H_\omega$ has AL for a.e. $\omega\in\tilde{\Omega}$, i.e. there exist $\{u_k(n,\omega)\}_{k\in\Z}$ which forms an orthonormal basis of $\ell^2(\Z^{d})$. Thus for $n>N_0$,
$$
|u_k(n,\omega)u_k(0,\omega)|\leq C_\omega e^{-(\ln\lambda+\frac{\e_0}{2})|n|}.
$$
Since $\sum_k|u_k(0,\omega)|^2=1$, there exists $k_0$ such that $u_{k_0}(0,\omega)\neq 0$, thus for $n>N_0$, we have
\begin{equation}\label{lower}
|u_{k_0}(n,\omega)|\leq \frac{C_\omega}{|u_{k_0}(0,\omega)|}e^{-(\ln\lambda+\frac{\e_0}{2})|n|}.
\end{equation}
This implies that
$$
\liminf_{n\rightarrow \infty} -\frac{\ln(|u_{k_0}(n,\omega)|^2+|u_{k_0}(n+1,\omega)|^2)}{2|n|}\geq \ln\lambda+\frac{\e_0}{2}.
$$
Note that for AMO, it is proved in \cite{J} that
$$
\lim_{n\rightarrow \infty} -\frac{\ln(|u_{k_0}(n,\omega)|^2+|u_{k_0}(n+1,\omega)|^2)}{2|n|}=\ln\lambda.
$$
This is a contradiction.
\qed \\


\noindent
\textbf{Proof of Theorem \ref{long}}.
Replace Proposition \ref{13}  by Proposition \ref{long2},  using similar argument as in Theorem \ref{subedl}, one can prove that  $L_{V, \lambda,\alpha,\theta}$ has EDL with $\gamma(L) \geq 2\pi(h-\e)$. We remark that in this case  we cann't obtain $\gamma(L) \geq 2\pi h$  and the largeness of $\lambda$ will depend on $\e$,  since there is no global almost reducibility result as in Lemma \ref{Global to local general}.
\qed \\

\noindent
\textbf{Proof of Theorem \ref{highamo}}. The proof is similar to that of Theorem \ref{long}. Note that $\{L_{\lambda,\alpha,\theta}\}_{\theta\in \T}$ is the duality of the following Schr\"odinger operator on $\ell^2(\Z)$:
\begin{equation}\label{highamodual}
(H_{\lambda^{-1} V,\alpha,\theta}u)_n=u_{n+1}+u_{n-1}+2\lambda^{-1}\sum\limits_{i=1}^d\cos2\pi(\theta_i+n\alpha_i)u_n,
\end{equation}
where $\theta=(\theta_1,\cdots,\theta_d)\in\T^d$, $\alpha=(\alpha_1,\cdots,\alpha_d)\in \T^d$.

One can directly see that $|2\lambda^{-1}\sum\limits_{i=1}^d\cos2\pi(\theta_i+n\alpha_i)|_h=\lambda^{-1}e^{2\pi h}$, thus let $h=\frac{1}{2\pi}(1-\frac{\e}{2})\ln\lambda$, $\tilde{h}=\frac{1}{2\pi}(1-\frac{3\e}{4})\ln\lambda$ and $\lambda_0$
sufficiently large such that
\begin{equation}\label{initial2}
|2\lambda_0^{-1}\sum\limits_{i=1}^d\cos2\pi(\theta_i+n\alpha_i)|_h=\lambda_0^{-\frac{\e}{2}}\leq \frac{D_0}{\|A_E\|^{C_0}}(\frac{\e}{4\pi})^{C_0\tau},
\end{equation}
where $C_0,D_0$ are the constants defined in  Proposition \ref{iteration} and $A_E=\begin{pmatrix}E&-1\\1&0\end{pmatrix}$. $\|A_E\|$ is uniformly bounded for $E$ in the spectrum, thus $\lambda_0$ only depends on $\alpha,d,\e$.

Note that \eqref{initial2} implies that Proposition \ref{reducibility} holds for the associated Schr\"odinger cocycle of the operator \eqref{highamodual} provided $\lambda>\lambda_0$. Then we use Propostion \ref{long2} and similar argument as in Theorem \ref{subedl} to prove that  $L_{\lambda,\alpha,\theta}$ has EDL with $\gamma(L) \geq2\pi(\tilde{h}-\frac{\e}{96})\geq (1-\e)\ln\lambda$.\qed

\section{Appendix}

The  following quantitative almost reducibility result is the basis of our proof. The result first appeared in \cite{LYZZ}, which is refined version of  \cite{houyou,ccyz}.

\begin{Proposition}\label{iteration}
Let $\alpha\in DC(\kappa',\tau)$, $\sigma>0$. Suppose that $A\in SL(2,\R)$, $f\in C^\omega_h(\T^d,\sl(2,\R))$. Then for any $h_+<h$, there exists numerical constant $C_0$,  and constant $D_0=D_0(\kappa',\tau,d)$ such that if

\begin{align}\label{initial}
\|f\|_h\leq \epsilon\leq \frac{D_0}{\|A\|^{C_0}}(\min\{1,\frac{1}{h}\}(h-h_+))^{C_0\tau},
\end{align}
then there exists $B\in C_{h_+}^\omega(2\T^d,SL(2,\R))$, $A_+\in SL(2,\R)$ and $f_+\in C_{h_+}^\omega(\T^d,\sl(2,\\\R))$ such that
$$
B^{-1}(\theta+\alpha)Ae^{f(\theta)}B(\theta)=A_+e^{f_+(\theta)}.
$$
More precisely, let $spec(A)=\{e^{2\pi i\xi},e^{-2\pi i\xi}\}$, $N=\frac{2}{h-h_+} | \ln \epsilon |$, then we can distinguish two cases:
\begin{itemize}
\item (Non-resonant case)   if for any $n\in \Z^{d}$ with $0< |n| \leq N$, we have
$$
\| 2\xi - <n,\alpha> \|_{\R/\Z} \geq \epsilon^{\frac{1}{15}},
$$
then
$$\| B-id\|_{h_+}\leq \epsilon^{\frac{1}{2}} , \quad   \|f_{+}\|_{h_+}\leq \epsilon^{2}.$$
Moreover, $\|A_+-A\|<2\epsilon$.
\item (Resonant case) if there exists $n_\ast$ with $0< |n_\ast| \leq N$ such that
$$
\| 2\xi- <n_\ast,\alpha> \|_{\R/\Z}< \epsilon^{\frac{1}{15}},
$$
then $B(\theta)=B'R_{\langle n_*,\theta\rangle}$ with estimates
$$ \|B\|_{h_+}\leq |\ln \epsilon |^{\tau}\epsilon^{-\frac{h_+}{h-h_+}},\ \  \|B'\|_{h_+}\leq | \ln \epsilon |^{\tau}, \ \ \|f_{+}\|_{h_+}< \epsilon e^{-h_+\epsilon^{-\frac{1}{18\tau}}}.$$
Moreover, $\deg B=n_*$,  the constant  $ A_+$  can be written as
 $$A_+= M^{-1}exp \left(
\begin{array}{ccc}
 i t^{+} &  \nu^{+}\\
\bar{ \nu}^{+} &  -i t^{+}
 \end{array}\right)M
 $$ with estimates  $|\nu^{+}|\leq e^{-2\pi|n_*|h}$, $|t^{+}|\leq \epsilon^{\frac{1}{16}}$.
\end{itemize}
\end{Proposition}

\begin{Remark}
Note that if $h<1$, Proposition \ref{iteration} was proved in \cite{ccyz,LYZZ}. If $h>1$, to ensure the inequality in the top of page 17 in \cite{LYZZ}
$$
|P|_{h_+}\leq (d-1)!(\mathcal{N}_2+\frac{1}{h-h_+})^d\cdot\tilde{\epsilon}e^{-2\pi(h-h_+)\mathcal{N}_2}\leq \epsilon e^{-h_+\epsilon^{-\frac{1}{16\tau}}}
$$
where $\mathcal{N}_2=2^{-\frac{1}{\tau}}\gamma^{\frac{1}{\tau}}\epsilon^{-\frac{1}{15\tau}}-N$ holds, one needs to replace
$$
\epsilon<\frac{D_0}{\|A\|^{C_0}}(h-h_+)^{C_0\tau}
$$
by
$$
\epsilon<\frac{D_0}{\|A\|^{C_0}}(\frac{h-h_+}{h})^{C_0\tau}.
$$
The rest estimates  are the same as in \cite{ccyz,LYZZ}.
\end{Remark}

\begin{Theorem}\label{positive reducibility}
Let $\alpha\in DC_d(\kappa',\tau)$, $h>\tilde{h}>0$, $\tau>d-1$, $\kappa'>0$, $\kappa>0$,  $R\in SL(2,\R)$.   Let   $A\in C_{h}^\omega(\T^d,SL(2,\R))$ with   $rot_f(\alpha,A(\theta))\in DC_\alpha(\kappa,\tau)$, where
$$
DC_\alpha(\kappa,\tau)=\{\phi\in\R^d|\|2\phi-m\alpha\|_{\R/\Z}\geq \frac{\kappa}{(|m|+1)^{\tau}}\}
$$
Then  there exist numerical constant $C_0$, constant $D_0=D_0(\kappa',\tau,d)$,  $\epsilon=\epsilon(\tau,\kappa',\kappa,h,\tilde{h},d,R)$, such that if
$$
\|A(\theta)-R\|_{h}\leq \epsilon  \leq  \frac{D_0\kappa^{4}}{\|A\|^{C_0}}\min\{1,\frac{1}{h}\}(h-\tilde{h})^{C_0\tau},
$$
then there exist $B\in C_{\tilde{h}}^\omega(\T^d,SL(2,\R))$,  $\tilde{A}\in SL(2,\R)$ such that
$$
B(\theta+\alpha)A(\theta)B(\theta)^{-1}=\tilde{A},
$$
with estimates $\|B-id\|_{\tilde{h}}\leq \|A(\theta)-R\|_{h}^{\frac{1}{2}}$,  $\|\tilde{A}-R\|\leq \|A(\theta)-R\|_{h}$ .
\end{Theorem}

This result was essentially proved by Dinaburg-Sinai \cite{ds}. Here we sketch the proof since we need  a bit more concrete estimates.
\begin{pf} Assume that  $\epsilon_0$ is small enough such that
\begin{equation}\label{f1}
\|A(\theta)-A\|_{h_0}\leq \epsilon_0:=\frac{D_0\kappa^{4}}{\|A\|^{C_0}}\min\{1,\frac{1}{h}\}(h-\tilde{h})^{C_0\tau},
\end{equation}
we can write $A(\theta)=Ae^{f(\theta)}$ with $\|f(\theta)\|_{h}\leq\epsilon_0$.

We prove by induction. Take $\epsilon_0$, $h$ and $\tilde{h}$ as above. Assume that we are at the $(j+1)^{th}$ KAM step, where we have $A_j\in SL(2,\R)$ with two eigenvalues $e^{\pm i\xi_j}$ and $f_j\in \mathfrak{B}_{h_j}$ satisfying $\|f_j\|_{h_j}\leq \epsilon_j$ for some $\epsilon_j\leq \epsilon_0^{2^j}$, then we define
$$
h_j-h_{j+1}=\frac{h-\tilde{h}}{4^{j+1}},\ \ N_j=\frac{2|\ln\epsilon_j|}{h_j-h_{j+1}}.
$$

By $(\ref{f1})$, it is easy to check that
$$
\epsilon_j\leq\frac{D_0\kappa^4}{\|A\|^{C_0}}\min\{1,\frac{1}{h_j}\}(h_j-h_{j+1})^{C_0\tau}.
$$
Thus for any $n\in\Z^d$ with $0<|n|<N_j$, we have
\begin{align*}
\|\xi_{j}-\langle n,\alpha\rangle\|_{\R/\Z}&\geq -|\xi_{j}-rot_f(\alpha,A_j+F_j)|+\|rot_f(\alpha,A_j+F_j)-\langle n,\alpha\rangle\|_{\R/\Z}\\
&\geq \frac{\kappa}{(|n|+1)^{\tau}}-\epsilon_{j}\geq \frac{\kappa}{(|N_j|+1)^{\tau}}-\epsilon_{j}\geq \epsilon_j^{\frac{1}{2}}.
\end{align*}

Let $\epsilon=\epsilon_j$, $h=h_j$, $h_+=h_{j+1}$ and $A=A_j$, then, by Proposition \ref{iteration}, we can construct
$$
B_j\in C^\omega_{h_{j+1}}(2\T^d,SL(2,\R)),\ \ A_{j+1}\in SL(2,\R),\ \ f_{j+1}\in \mathfrak{B}_{h_{j+1}},
$$
such that
$$
B_j(\theta+\alpha)A_je^{f_j(\theta)}B_j(\theta)^{-1}=A_{j+1}e^{f_{j+1}(\theta)}.
$$
with
$$
\|A_{j+1}-A_j\|\leq \epsilon_j^{\frac{1}{2}},\ \ \|B_j-id\|_{h_{j+1}}\leq 2\epsilon_j^{\frac{1}{2}}, \ \ \|f_{j+1}\|_{h_{j+1}}\leq \epsilon_{j+1}=\epsilon_j^2.
$$
Let $B(\theta)=\prod_jB_j(\theta)$, then $B\in C_{\tilde{h}}^\omega(\T^d,SL(2,\R))$, and
$$
B(\theta+\alpha)A(\theta)B(\theta)^{-1}=\tilde{A},
$$
with estimates $\|B-id\|_{\tilde{h}}\leq \|A(\theta)-A\|_{h}^{\frac{1}{2}}$.
\end{pf}

\section*{Acknowledgements}
The authors want to thank S.Jitomirskaya for useful discussions.
J. You  was partially supported by NSFC grant (11871286) and
Nankai Zhide Foundation.  Q. Zhou was partially supported by  NSFC grant (11671192,11771077)  and
Nankai Zhide Foundation.

\end{document}